\documentclass[nohyperref]{article}

\usepackage{microtype}
\usepackage{graphicx}
\usepackage{subfigure}
\usepackage{booktabs} 
\usepackage{multirow}
\usepackage{diagbox}

\usepackage{hyperref}

\usepackage[accepted]{icml2023}

\usepackage{amsmath}
\usepackage{amssymb}
\usepackage{mathtools}
\usepackage{amsthm}

\usepackage[capitalize,noabbrev]{cleveref}

\theoremstyle{plain}

\theoremstyle{definition}

\theoremstyle{remark}

\usepackage[textsize=tiny]{todonotes}

\icmltitlerunning{General Covariance Data Augmentation}

\begin{document}

\twocolumn[
\icmltitle{General Covariance Data Augmentation for Neural PDE Solvers}



\icmlsetsymbol{equal}{*}

\begin{icmlauthorlist}
\icmlauthor{Vladimir Fanaskov}{Sk}
\icmlauthor{Tianchi Yu}{Sk}
\icmlauthor{Alexander Rudikov}{Sk,INM}
\icmlauthor{Ivan Oseledets}{Sk,AIRI}
\end{icmlauthorlist}

\icmlaffiliation{Sk}{Skoltech, Center for Artificial Intelligence Technology}
\icmlaffiliation{AIRI}{Artificial Intelligence
Research Institute}
\icmlaffiliation{INM}{Marchuk Institute of Numerical Mathematics of the Russian Academy of Sciences}

\icmlcorrespondingauthor{Vladimir Fanaskov}{V.Fanaskov@skoltech.ru}

\icmlkeywords{Machine Learning, ICML}

\vskip 0.3in
]



\printAffiliationsAndNotice{} 

\begin{abstract}
    The growing body of research shows how to replace classical partial differential equation (PDE) integrators with neural networks. The popular strategy is to generate the input-output pairs with a PDE solver, train the neural network in the regression setting, and use the trained model as a cheap surrogate for the solver. The bottleneck in this scheme is the number of expensive queries of a PDE solver needed to generate the dataset. To alleviate the problem, we propose a computationally cheap augmentation strategy based on general covariance and simple random coordinate transformations. Our approach relies on the fact that physical laws are independent of the coordinate choice, so the change in the coordinate system preserves the type of a parametric PDE and only changes PDE's data (e.g., initial conditions, diffusion coefficient). For tried neural networks and partial differential equations, proposed augmentation improves test error by 23\% on average. The worst observed result is a 17\% increase in test error for multilayer perceptron, and the best case is a 80\% decrease for dilated residual network.
\end{abstract}

\section{Introduction}
Machine learning is increasingly used to solve partial differential equations (PDEs). The especially fruitful idea is to learn a computationally cheap but sufficiently accurate surrogate for the classical solver \cite{hennigh2017lat}, \cite{li2020fourier}, \cite{tripura2022wavelet}, \cite{lu2021learning}, \cite{stachenfeld2021learned}. The most reliable training strategy is to generate input-output pairs with a classical solver and fit a neural network of choice with a standard $L_2$ loss (regression setting).

An alternative we do not consider here is to resort to a so-called physics-informed setting when the loss is a $L_2$ norm of PDE residual evaluated at certain points \cite{wang2021learning}, \cite{li2021physics}. This way, one avoids data generation by a classical solver. Arguably, it is currently recognized that the whole process is inefficient \cite{lu2021deepxde}, \cite{karnakov2022solving}.

In the regression setting the size of the generated dataset is usually limited owing to the restrictions on the computation budget. Deep learning is data-hungry, so ways to cheaply increase the number of data points available for training are highly desirable. Numerous augmentation techniques serve this purpose in classical machine learning \cite{shorten2019survey}, \cite{wen2020time}. For scientific machine learning, literature on augmentation is scarce \cite{brandstetter2022lie}, \cite{li2022physics}. In this note, we contribute a new way to augment datasets for neural PDE solvers.

The central idea behind our approach is \textbf{the principle of general covariance}. General covariance states that physical phenomena do not depend on the choice of a coordinate system \cite{post1997formal}, \cite{emam2021covariant}. Mathematically, the covariance means the physical fields are geometric objects (tensors) with particular transformation laws under the change of coordinates \cite{MR1441306}, \cite{MR3675718}. In exceptional cases, these transformation laws leave governing equations invariant (symmetry transformation), but in most cases, it is only the form of the equations that persists. More specifically, for parametric partial differential equations, suitably chosen coordinate transformation induces the change of problem data (e.g., permeability field, convection coefficients, source term, initial or boundary conditions, e.t.c.). We use this fact to build a computationally cheap and broadly applicable augmentation strategy based on simple random coordinate transformations. To evaluate the efficiency of our approach, we perform empirical tests on the two-way wave, convection-diffusion, and stationary diffusion equations using several variants of Fourier Neural Operator (FNO) \cite{li2020fourier}, Deep Operator Network (DeepONet) \cite{lu2021learning}, Multilayer Perceptron (MLP) \cite{haykin1994neural}, Dilated Residual Network (DilResNet) \cite{yu2015multi}, \cite{stachenfeld2021learned} and U-Net \cite{ronneberger2015u}. Both for one-dimensional and two-dimensional PDEs proposed augmentation technique improves test error by 23\% on average and up to 80\% in the most favorable cases.

\textbf{Contributions:}
\begin{enumerate}
	\item Easily extendable, architecture-agnostic augmentation procedure based on general covariance.
	\item Cheap algebraic random grids in $R^{n}$ based on cumulative distribution function and transfinite interpolation.
	\item Comprehensive set of experiments showing that augmentation helps to improve test error for different architectures and parametric families of PDEs.
\end{enumerate}

Code and datasets are available on \url{https://github.com/VLSF/augmentation}.

\label{section:introduction}
\begin{figure*}[t]
\vskip 0.2in
\begin{center}
\centerline{\includegraphics[width=\textwidth]{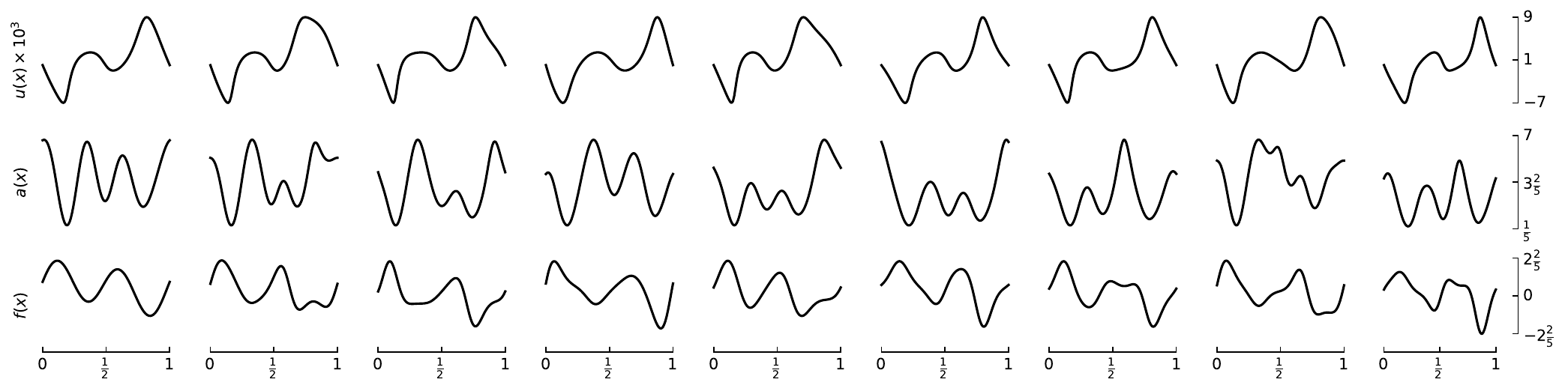}}
\caption{Example of data augmentation for elliptic equation \eqref{eq:elliptic_1D}. The first column on the left contains features and target that solves \eqref{eq:elliptic_1D}. All other columns are obtained from the first one with transformation \eqref{eq:elliptic_augmentation}. Coordinate transformations are generated according to \eqref{eq:coordinate_transformation} with parameters $N=3$, $\beta=1$ and coefficients $c_k, d_k,\,k=1,2,3$ sampled from standard normal distribution.}
\label{fig:coordinate_transform}
\end{center}
\vskip -0.2in
\end{figure*}

\section{Basic augmentation example}
\label{section:Basic augmentation example}
Before dwelling upon technical details, we provide a simple example of our approach for a parametric boundary-value problem
\begin{equation}
    \label{eq:elliptic_1D}
    \begin{split}
    &\frac{d}{dx}\left(a(x) \frac{d}{dx} u(x)\right) = f(x),\\
    &x\in[0, 1],\,u(0) = u(1) = 0.
    \end{split}
\end{equation}

Suppose that $a$ and $f$ are chosen reasonably, so the unique solution exists.\footnote{The formal statement on existence is available in \cite{evans2010partial}, but it is largely irrelevant to our discussion.} The usual way to approximate this solution is to use finite-element discretization
\begin{equation}
    \label{eq:FEM}
    u(x) = \sum_{i=1}^{N} \phi_i(x) u_i,
\end{equation}
where $\phi_i(x)$ are piecewise linear functions such that $\phi_i(x_j) = \delta_{ij}$ for $x_j = j/(N+1),\,j=1,\dots,N$, i.e., the hat functions. After that, the differential equation \eqref{eq:elliptic_1D} in a weak form is equivalent to $N\times N$ linear system, and the solution is straightforward. The obtained solution is known not only on the uniform grid $\mathcal{G} = \left\{j/(N+1),\,j=0,\dots, N+1\right\}$ but everywhere in the domain thanks to the closed-form representation \eqref{eq:FEM}.

Using described procedure one can generate a dataset with features $F_i = \left(a_i(\mathcal{G}), f_i(\mathcal{G})\right)$ and targets $T_i = \left(u_i(\mathcal{G})\right)$, $i=1,\dots, N_{\text{samples}}$. As a rule, features, $a_i$ and $f_i$ in our case, are samples from some distribution \cite{kovachki2021neural} or typical inputs needed for a particular application, e.g., \cite{pathak2022fourcastnet}.

Our main observation is that when the PDE is known, it is possible to extract more information from each obtained solution using coordinate transformation. Suppose $y(\xi)$ is analytic strictly monotonic function from $[0, 1]$ to $[0, 1]$ such that $y(0) = 0$, $y(1) = 1$. We use $x \equiv y(\xi)$ as coordinate transformation and rewrite \eqref{eq:elliptic_1D} in coordinates $\xi$ as follows
\begin{equation}
    \label{eq:elliptic_1D_transformed}
    \begin{split}
    &\frac{d\xi}{d y}\frac{d}{d\xi}\left(a(y(\xi)) \frac{d\xi}{d y}\frac{d}{d\xi} u(y(\xi))\right) = f(y(\xi)),\\
    &\xi\in[0, 1],\,u(y(0)) = a,\,u(y(1)) = b.
    \end{split}
\end{equation}
As we see, transformed equation \eqref{eq:elliptic_1D_transformed} has the same parametric form as the original one \eqref{eq:elliptic_1D}. As a consequence, if a triple of functions $a(x), u(x), f(x)$ solve \eqref{eq:elliptic_1D}, the triple of modified functions $a(y(x))  \frac{dx}{dy}, u(y(x)), f(y(x)) \frac{dy}{dx}$ also solve the same equation \eqref{eq:elliptic_1D}, where we rename variable $\xi$ in \eqref{eq:elliptic_1D_transformed} to $x$. So we can generate novel solutions from the old ones using smooth coordinate transformations and interpolation.

To complete a description of the augmentation, we need to explain how to generate smooth coordinate transformations. Since any strictly monotonic positive function that maps $[0, 1]$ constitutes a valid coordinate transformation, we propose to use cumulative distribution functions with strictly positive probability density. It is easy to come up with many parametric families of probability densities. For example, we can use trigonometric series and define
\begin{equation}
    \begin{split}
    &p(x) = 1 + \sum_{k=1}^{N}\frac{\left(c_k \cos(2\pi k x) + d_k\sin(2\pi k x)\right)}{c_0},\\
    & c_0 = \sum_{k=1}^{N}\left(\left|c_k\right| + \left|d_k\right|\right) + \beta,\,\beta>0.
    \end{split}
\end{equation}
After integration, we obtain a cumulative distribution function that serves as a coordinate transformation
\begin{equation}
    \label{eq:coordinate_transformation}
    y(x) = x + \sum_{k=1}^{N}\frac{\left(c_k \sin(2\pi k x) + d_k(1-\cos(2\pi k x))\right)}{2\pi k c_0}.
\end{equation}
The whole augmentation procedure for elliptic equation \eqref{eq:elliptic_1D} can be compactly written as
\begin{equation}
    \label{eq:elliptic_augmentation}
    \underset{\text{solve \eqref{eq:elliptic_1D}}}{a(x), u(x), f(x)} \longrightarrow \underset{y(x) \text{ from \eqref{eq:coordinate_transformation}}}{a(y(x))\big/\frac{dy}{dx}, u(y(x)), f(y(x)) \frac{dy}{dx}}.
\end{equation}
\cref{fig:coordinate_transform} illustrates the proposed approach for elliptic equation \eqref{eq:elliptic_1D} and a particular set of transformations \eqref{eq:coordinate_transformation}.

To summarize, our augmentation approach consists of three stages:
\begin{enumerate}
    \item Generate a sufficiently smooth coordinate transformation $y(\xi)$.
    \item Interpolate features and targets on discrete grid $y(\xi_j),\,\xi_j=j/(N+1),\,j=0,\dots,N+1$.
    \item Adjust interpolated features and targets according to the transformations law for PDE evaluated in new coordinates $y(\xi)$.
\end{enumerate}
This procedure can be applied for as many coordinate transformations $y(\xi)$ as needed and requires only cheap interpolation, so the overall cost is $O(N)$ for each sample, where $N$ is the number of grid points.

In the \cref{section:Augmentation by General Covariance}, we show how to generalize results illustrated here for other partial differential equations and higher dimensions.

\label{section:augmentation example}
\section{Augmentation by General Covariance}
\label{section:Augmentation by General Covariance}
In \cref{section:Basic augmentation example}, we explained that two principal components of the proposed augmentation approach are grid generation and transformation law for PDE in question. Here we show how to extend results from \cref{section:Basic augmentation example} to a more general setting.

Everywhere in this section, Einstein's summation notation is used, e.g., $a_{\alpha}b^{\alpha} \equiv \sum_{\alpha} a_{\alpha}b^{\alpha}$.

\subsection{How to construct coordinate transformations in the general case}
\label{section:Augmentation by General Covariance.subsection:Grid generation}
We define coordinate transformations in $D\geq1$ as one-to-one analytic mapping
\begin{equation}
    \label{eq:coordinate transformation}
    \boldsymbol{x}(\boldsymbol{\xi}):[0, 1]^{D}\longrightarrow[0, 1]^{D}
\end{equation}

In \cref{section:Basic augmentation example}, we outlined a particular scheme to construct families of coordinate transformations in $D=1$. The general algorithm is as follows
\begin{enumerate}
    \item Select a family of basis functions $\phi_j(\xi)$ defined on $[0, 1]$ that are easy to integrate (e.g., the indefinite integral is known).
    \item Find suitable shift and scale for a series $s\left(\sum_{j}\phi_j(\xi)c_{j}\right) + c_0$ to be a valid probability density function $p(\xi)$ for all $c_j$.
    \item Use cumulative distribution function (indefinite integral of $p(\xi)$) as a coordinate transformation.
\end{enumerate}

When $D=1$ mapping is available, it is possible to lift it to $D>1$ by transfinite interpolation \cite{gordon1973construction}. For example, for $D=2$ the transformation becomes
\begin{equation}
    \label{eq:2D mapping}
    \begin{split}
        &x^{1}(\xi^{1}, \xi^{2}) = y_1(\xi^{1})(1 - \xi^{2}) + y_2(\xi^{1})\xi^{2},\\
        &x^{2}(\xi^{1}, \xi^{2}) = y_3(\xi^{2})(1 - \xi^{1}) + y_4(\xi^{2})\xi^{1},
    \end{split}
\end{equation}
where $y_i,\,i=1,\dots,4$ are $D=1$ mappings, e.g., given in \eqref{eq:coordinate_transformation}. The extension of \eqref{eq:2D mapping} to higher dimensions is straightforward.

Note that \eqref{eq:2D mapping} has a ``low-rank'' structure that decreases the diversity of possible grids. The issue can be alleviated with Hermite transfinite interpolation or with more general blending functions \cite{MR3675718}. Other more computationally involved remedies are variational and elliptic (Laplace-Beltrami) grid generators \cite{steinberg1986variational}, \cite{spekreijse1995elliptic}.

\subsection{Linear PDEs under coordinate transformations}
Here we remind how the most widely used differential operators change under the transformation \eqref{eq:coordinate transformation}.

The results we present in this section are standard \cite{MR3675718}, \cite{MR1441306}, \cite{MR1247707}. For convenience, the proofs are also available in \cref{appendix:Coordinate transformations}.

For convenience, we define the Jacobi matrix, and its determinant
\begin{equation}
    \mathcal{J}_{i\alpha} \equiv \frac{\partial x^{i}}{\partial \xi^{\alpha}},\,J \equiv \det \mathcal{J}.
\end{equation}
Note that for the mapping \eqref{eq:coordinate_transformation} determinant $J$ vanishes nowhere in the domain due to the strict monotonicity of the mapping.

It is not hard to show that for arbitrary space-dependent fields $c^{j}$, $\phi$, $a^{kj}$, $k, j\in 1,\dots, D$, the following transformation laws hold
\begin{equation}
    \label{eq:differential operators}
    \begin{split}
        &c^{j}\frac{\partial\phi}{\partial x^{j}} = c^{j}\frac{\partial \xi^{\alpha}}{\partial x^{j}}\frac{\partial\phi}{\partial \xi^{\alpha}},\\
        &a^{kj} \frac{\partial^2\phi}{\partial x^{j}\partial x^{k}} = a^{kj} \frac{\partial \xi^{\beta}}{\partial x^{j}} \frac{\partial \xi^{\gamma}}{\partial x^{k}}\frac{\partial^2 \phi}{\partial \xi^{\beta}\partial \xi ^{\gamma}} + a^{kj}\frac{\partial^2 \xi^{\gamma}}{\partial x^{k}\partial x^{j}} \frac{\partial \phi}{\partial \xi^{\gamma}},\\
        &\frac{\partial}{\partial x^{\alpha}}\left(c^{\alpha} \phi\right) = \frac{1}{J}\frac{\partial}{\partial \xi^{k}}\left(Jc^{\alpha}\frac{\partial \xi^{k}}{\partial x^{\alpha}}\phi\right),\\
        &\frac{\partial}{\partial x^{k}}\left(a^{kj}\frac{\partial\phi}{\partial x^{j}}\right) = \frac{1}{J}\frac{\partial}{\partial \xi^{k}}\left(J\left( a^{\alpha j}\frac{\partial \xi ^{k}}{\partial x^{\alpha}}  \frac{\partial \xi^{\beta}}{\partial x^{j}} \right)\frac{\partial\phi}{\partial \xi^{\beta}}\right).\\
    \end{split}
\end{equation}
Results \eqref{eq:differential operators} allow deriving transformation laws for many practically-relevant PDEs. We are interested in the following ones:
\begin{enumerate}
    \item Stationary diffusion equation
    \begin{equation}
        \label{eq:stationary diffusion}
        \begin{split}
            &\frac{\partial}{\partial x^{k}}\left(a^{kj}(\boldsymbol{x})\frac{\partial}{\partial x^{j}}u(\boldsymbol{x})\right) = f(\boldsymbol{x})\\
            &\boldsymbol{x}\in\Gamma \equiv [0, 1]^{D},\,\left.u(\boldsymbol{x})\right|_{x\in\partial\Gamma} = 0,
        \end{split}
    \end{equation}
    where $\partial\Gamma$ is a boundary of the unit hypercube $\Gamma$, and $\boldsymbol{x}\in\mathbb{R}^{D}$.
    \item Convection-diffusion equation
    \begin{equation}
        \label{eq:convection diffusion}
        \begin{split}
            \frac{\partial}{\partial t}\phi(\boldsymbol{x}, t) + &\frac{\partial}{\partial x^{i}} \left(v^{i}(\boldsymbol{x})\phi(\boldsymbol{x}, t)\right) = \\ &\frac{\partial }{\partial x^{k}}\left(a^{kj}(\boldsymbol{x}) \frac{\partial }{\partial x^{j}} \phi(\boldsymbol{x}, t)\right),\\
            \boldsymbol{x}\in\Gamma \equiv [0, &1]^{D},\,\left.\phi(\boldsymbol{x}, t)\right|_{x\in\partial\Gamma} = 0,\,\phi(\boldsymbol{x}, 0) = f(\boldsymbol{x}).
        \end{split}
    \end{equation}
    \item Two-way wave equation
    \begin{equation}
        \label{eq:wave equation}
        \begin{split}
            \frac{\partial^2 \rho(\boldsymbol{x}, t)}{\partial t^2} + v^{i}(\boldsymbol{x})\frac{\partial \rho(\boldsymbol{x}, t)}{\partial x^{i}} &= c^{ij}(\boldsymbol{x}) \frac{\partial^2 \rho(\boldsymbol{x}, t)}{\partial x^{i} \partial x^{j}} \\&+ e(\boldsymbol{x})\rho(\boldsymbol{x}, t),\\
            \boldsymbol{x}\in\Gamma \equiv [0, 1]^{D},\,\rho(\boldsymbol{x}, t\left.)\right|_{x\in\partial\Gamma} &= 0,\,\rho(\boldsymbol{x}, 0) = f(\boldsymbol{x}).
        \end{split}
    \end{equation}
\end{enumerate}
For all of the equations above, the transformed form easily follows from \eqref{eq:differential operators}. However, for wave and convection-diffusion equations, additional steps are required to ensure that the transformed equation has the same parametric form as the original one. \cref{table:PDEs transforms} contains the results for selected PDEs.

\begin{table*}[t]
    \centering
    \caption{Transformation of PDEs parameters under the change of coordinates $\boldsymbol{x}\rightarrow\boldsymbol{x}(\boldsymbol{\xi})$.}
    \vskip 0.15in
    \begin{tabular}{@{}lcc}
        \toprule
        Equation & Fields & Transformed fields\\ \cmidrule(r){1-3}
        \multirow{3}{*}{Stationary diffusion \eqref{eq:stationary diffusion}} & $u(\boldsymbol{x})$ & $u(\boldsymbol{x}(\boldsymbol{\xi}))$ \\ \addlinespace
        & $a^{k\beta}(\boldsymbol{x})$ & $Ja^{\alpha j}(\boldsymbol{x}(\boldsymbol{\xi}))\frac{\partial \xi ^{k}}{\partial x^{\alpha}} \frac{\partial \xi^{\beta}}{\partial x^{j}}$ \\ \addlinespace
        & $f(\boldsymbol{x})$ & $Jf(\boldsymbol{x}(\boldsymbol{\xi}))$ \\ \cmidrule(r){2-3}
        \multirow{4}{*}{Convection-diffusion \eqref{eq:convection diffusion}} & $\phi(\boldsymbol{x}, t)$ & $J \phi(\boldsymbol{x}(\boldsymbol{\xi}), t)$\\ \addlinespace
         & $a^{k\beta}(\boldsymbol{x})$ & $a^{\alpha j}(\boldsymbol{x}(\boldsymbol{\xi}))\frac{\partial \xi ^{k}}{\partial x^{\alpha}} \frac{\partial \xi^{\beta}}{\partial x^{j}}$\\ \addlinespace
         & $v^{i}(\boldsymbol{x})$ & $v^{k}(\boldsymbol{x}(\boldsymbol{\xi}))\frac{\partial \xi^{i}}{\partial x^{k}} + a^{\alpha j}(\boldsymbol{x}(\boldsymbol{\xi}))\frac{\partial \xi ^{k}}{\partial x^{\alpha}}  \frac{\partial \xi^{\beta}}{\partial x^{j}} \frac{\partial \xi^{\gamma}}{\partial x^{\rho}}\frac{\partial^2 x^{\rho}}{\partial \xi^{\gamma}\partial \xi^{\beta}}$\\ \addlinespace
         & $f(\boldsymbol{x})$ & $Jf(\boldsymbol{x}(\boldsymbol{\xi}))$\\ \cmidrule(r){2-3}
        \multirow{4}{*}{Wave \eqref{eq:wave equation}} & $\rho(\boldsymbol{x}, t)$ & $\rho(\boldsymbol{x}(\boldsymbol{\xi}), t)$\\ \addlinespace
         & $c^{\gamma\beta}(\boldsymbol{x})$ & $c^{kj}(\boldsymbol{x}(\boldsymbol{\xi})) \frac{\partial \xi^{\gamma}}{\partial x^{k}} \frac{\partial \xi^{\beta}}{\partial x^{j}}$\\ \addlinespace
         & $v^{\alpha}(\boldsymbol{x})$ & $v^{i}(\boldsymbol{x}(\boldsymbol{\xi}))\frac{\partial \xi^{\alpha}}{\partial x^{i}} - c^{kj}(\boldsymbol{x}(\boldsymbol{\xi}))\frac{\partial^2 \xi^{\alpha}}{\partial x^{k}\partial x^{j}} $\\ \addlinespace
         & $f(\boldsymbol{x})$ & $f(\boldsymbol{x}(\boldsymbol{\xi}))$\\
         & $e(\boldsymbol{x})$ & $e(\boldsymbol{x}(\boldsymbol{\xi}))$\\
        \bottomrule
    \end{tabular}
    \label{table:PDEs transforms}
    \vskip -0.1in
\end{table*}

Results summarized in \cref{table:PDEs transforms} along with the coordinate transformations described in \cref{section:Augmentation by General Covariance.subsection:Grid generation} are sufficient to perform general covariance augmentation for equations \eqref{eq:stationary diffusion}, \eqref{eq:convection diffusion}, \eqref{eq:wave equation}.

\subsection{Navier-Stokes equation}
To show that our scheme applies to nonlinear PDEs and more general boundary conditions, we consider lid-driven cavity flow of incompressible fluid with deformed cavities. The system of equations in the physical space reads
\begin{equation}
    \label{eq:cavity flow}
    \begin{split}
	&\frac{\partial v^{i}}{\partial t} = \frac{\partial}{\partial x^{k}}\left(-v^{k} v^{i} - p + \nu \frac{\partial v^{i}}{\partial x^k}\right),\,\frac{\partial v^{k}}{\partial x^{k}} =0,\\
    &\left.p(t,x)\right|_{x\in L} = 0,\,\left.\frac{\partial p(t, x)}{\partial x^{i}}n^{i}\right|_{x \in \Gamma} = 0,\\
	&\left.v^{1}(t, x)\right|_{x\in L} = 1,\,\left.v^{2}(t, x)\right|_{x\in C} = \left.v^{1}(t, x)\right|_{x\in \Gamma} = 0,
    \end{split}
\end{equation}
where $x$ belong to the interior of the curve $C(x^{1}, x^{2}) = L(x^{1}, x^{2})\cup \Gamma (x^{1}, x^{2})$, $L(x^{1}, x^{2})$ represents the lid and $\Gamma(x^{1}, x^{2})$ --- the rest of the cavity's boundary, $t\in[0, T]$.

For \cref{eq:cavity flow} we explicitly specify the form of cavity $x^{i}(\xi^{1}, \xi^{2}),\,i=1,2$ using curvilinear coordinates $\xi^{1}, \xi^{2}\in\left[0, 1\right]^2$ (see \cref{section:Experiments} for the description of the cavities used). When solution at $t=T$ is obtained, the only parameter of the PDE is the geometry itself which is fully specified by $x^{i}(\xi^{1}, \xi^{2}),\,i=1,2$.

In these circumstances general covariance augmentation simplifies. Namely, we use random coordinate transformations \cref{eq:2D mapping} to form additional mapping $\xi^{i}(\widetilde{\xi}^{i})$ and reinterpolate obtained solution $v^{i}(\xi^{1}, \xi^{2}, T)$ on the grid $\widetilde{\xi}^{i},\,i=1,2$.

Although transformation low of the equation Navier-Stokes equation \cref{eq:cavity flow} is not directly used for augmentation, we still need it to obtain a solution in the computational domain $\xi^{i},\,i=1,2$ (see \cref{appendix:Details on lid-driven cavity flow} for details).

In \cref{section:Experiments}, we show the empirical performance of the proposed augmentation scheme for both linear and nonlinear PDEs.



\section{Experiments}
\label{section:Experiments}
Here we present an empirical evaluation of augmentation by general covariance. For that purpose, we design several experiments in $D=1$ and $D=2$. We start with a description of the shared experiments' setup.

\subsection{Setup}
\subsubsection{Neural networks}
As a rule, neural PDE solvers are either Neural Operators or classical architectures used for image processing.\footnote{There are also hybrid methods (e.g., \cite{bar2019learning}), but we do not consider them here.} Since our approach is architecture-agnostic, we include results for both types of neural networks.

On the side of neural operators, we include original versions of DeepONet \cite{lu2021learning} and FNO \cite{li2020fourier}, implemented in \url{https://github.com/lu-group/deeponet-fno} and \url{https://github.com/neural-operator/fourier_neural_operator}, respectively. Besides, for the $D=1$ case, we also implemented an FNO-like operator dubbed rFNO with FFT replaced by pure real transform based on a complete trigonometric family and trapezoidal rule. In addition, for $D=2$ we also use Spectral Neural Operator described in \cite{fanaskov2022spectral}. Roughly speaking, the architecture has the same structure as FNO but with Discrete Cosine Transform in place of FFT.

Classical machine-learning architectures include DilResNet \cite{yu2015multi}, \cite{stachenfeld2021learned}, U-Net \cite{ronneberger2015u}, and MLP \cite{haykin1994neural}.

A detailed description of neural networks is available in \cref{appendix:Architectures and training details}.

\subsubsection{Partial differential equations and datasets}
We evaluate augmentation on stationary diffusion \eqref{eq:stationary diffusion}, convection-diffusion \eqref{eq:convection diffusion}, and wave \eqref{eq:wave equation} equations.

To produce PDE data, for linear PDEs, we sampled all needed functions from random trigonometric series
\begin{equation}
    f(x) = \sum_{k=0}^{N-1}\left(c_k \cos(2\pi k x) + s_k \sin(2\pi k x)\right),
\end{equation}
with $c_k$ sampled from the standard normal distribution and scaled/shifted appropriately to ensure needed boundary conditions or make $f(x)$ uniformly positive. For $D=2$, the procedure is the same, but a direct product of one-dimensional bases is used.

Afterward, equations with randomly generated data are discretized either with finite-difference or finite-element methods.

For $D=1$, we generated one dataset per equation and an additional dataset for the wave equation. Results for the extra wave dataset are available in \cref{appendix:Datasets}. For $D=2$, we produced two datasets per equation that differ by complexity (more rough targets or more diverse feature-target pairs).

Also, for the purposes explained later, we use two distinct elliptic datasets in $D=2$. In the first one, named ``Elliptic alpha'', diffusion coefficients $a^{i\beta}$ form a symmetric positive definite matrix for each point of the domain. In the second one, named ``Elliptic beta'', the matrix $a^{i\beta}$ is the identity matrix multiplied by a single uniformly positive diffusion coefficient.

For Navier-Stokes equation \cref{eq:cavity flow} we generate cavities defined with a transfinite interpolation \cref{eq:2D mapping} from randomly generated boundary curves
\begin{equation}
    \begin{split}
        &x^{1}(\xi^{1}, 0) = \xi^{1},\,x^{1}(\xi^{1}, 1) = \xi^{1},\\
        &x^{2}(0, \xi^{2}) = \xi^{2},\,x^{2}(1, \xi^{2}) = \xi^{2},\,x^{2}(\xi^{1}, 1) = 1,\\
        &x^{1}(0, \xi^{2}) = \sum_{i=1}^{m}\sin(\pi \xi^{2} (k+1) c^k \big/(10(k+1)^2),\\
        &x^{1}(1, \xi^{2}) = 1 + \xi^{2}(1-\xi^{2})\alpha/2,\\
        &x^{2}(\xi^{1}, 0) = \sum_{i=1}^{m}\sin(\pi \xi^{1} (k+1) d^k \big/(10(k+1)^2),\\
    \end{split}
\end{equation}
where $d^{k}, c^{k}, \alpha$ are sampled from standard normal distribution and $\xi^{1}, \xi^{2}\in[0, 1]^2$. We also use $T=10^{-4}$ and $\nu=10^{-2}$.

More details on the dataset generation process are available in \cref{appendix:Datasets}. Links to the datasets are available in the repository \url{https://github.com/VLSF/augmentation}.

\subsubsection{Coordinate transformations}
To generate coordinate transformation for $D=1$, we use a cumulative distribution function constructed from unnormalize probability density
\begin{equation}
    p(x) = \beta + \sum_{i=1}^{N}\left(c_i \cos(2n\pi x + p_i)\right),
\end{equation}
where $c_i$ and $p_i$ are samples from the standard normal distribution, and $N = 5$, $\beta = 1$.

For $D=2$ we use four $D=1$ coordinate transformations constructed with \eqref{eq:coordinate_transformation}, where coefficients $i=1,\dots, 5$ are from standards normal distribution, $c_0=1$, $N = 6$, and $\beta = 10^{-5}$. To obtain $D=2$ mappings from the four unidimensional, we apply transformation using transfinite interpolation \eqref{eq:2D mapping}.

\subsubsection{Metrics}
As a main measure of performance, we use average relative $L_2$ test error
\begin{equation}
    \label{eq:L2 error}
    E_{\text{test}} = \frac{1}{N} \sum_{i=1}^{N} \frac{\left\|\mathcal{N}(f_i) - t_i\right\|_2}{\left\|t_i\right\|_2},
\end{equation}
where $\mathcal{N}$ is a neural network, $f_i$ and $t_i$ are features and targets from the test set.

To evaluate the impact of augmentation, we consider the relative gain
\begin{equation}
    \label{eq:relative gain}
    g = \left(1 - \frac{E_{\text{test}}^{\text{aug}}}{E_{\text{test}}}\right)\times 100\%,
\end{equation}
where $E_{\text{test}}^{\text{aug}}$, and $E_{\text{test}}>0$ are relative errors of neural network trained with and without augmentation respectively. Since $E_{\text{test}}^{\text{aug}} = (1 - g)E_{\text{test}}$, the larger $g$ the better.

\begin{table}[b]
    \centering
    \caption{Relative gain for averaging with respect to equations and networks. Augmentation factor $m$ means that additional $m N_{\text{train}}$ augmented samples are appended to the training dataset.}
    \vskip 0.15in
    \begin{tabular}{@{}ccccc@{}}
        \toprule
        $m\backslash N_{\text{train}}$    & $500$ & $1000$ & $1500$ & $2000$ \\ \cmidrule{1-5}
        $1$ & $11\%$ & $18\%$ & $21\%$ & $21\%$ \\
        $2$ & $16\%$ & $25\%$ & $27\%$ & $31\%$ \\
        $3$ & $15\%$ & $23\%$ & $28\%$ & $28\%$ \\
        $4$ & $19\%$ & $25\%$ & $28\%$ & $30\%$ \\
        \bottomrule
    \end{tabular}
    \label{table:augmentation-test}
    \vskip -0.1in
\end{table}
\begin{table}[b]
    \centering
    \caption{Relative gain for averaging wrt different training scenarios.}
    \vskip 0.15in
    \resizebox{\columnwidth}{!}{\begin{tabular}{@{}lccccc@{}}
        \toprule
         & DeepONet & FNO & DilResNet & rFNO & MLP\\ \cmidrule{2-6}
        elliptic & $16\%$ & $22\%$ & $28\%$ & $19\%$ & $8\%$\\
        conv-dif & $32\%$ & $36\%$ & $22\%$ & $24\%$ & $28\%$\\
        wave & $14\%$ & $36\%$ & $22\%$ & $18\%$ & $16\%$\\
        \bottomrule
    \end{tabular}}
    \label{table:equation vs network}
    \vskip -0.1in
\end{table}

\subsection{Sensitivity to grid distortions}

\begin{figure}[tb]
    \begin{center}
    \centerline{\includegraphics[width=\columnwidth]{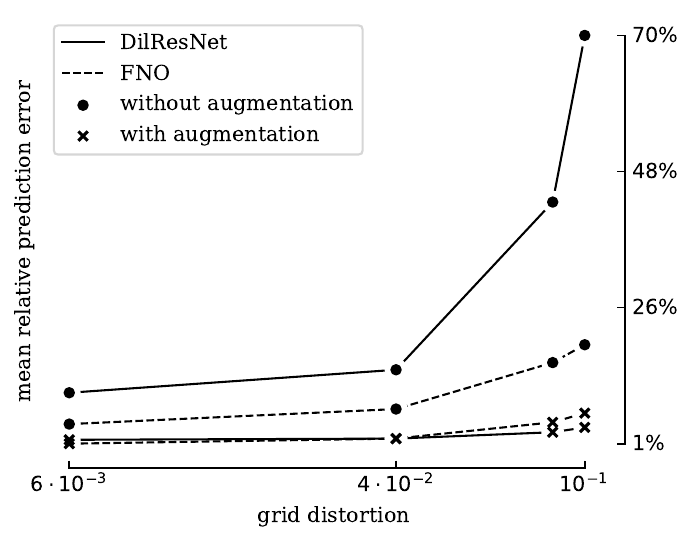}}
    \caption{Sensitivity to grid distortion for DilResNet and FNO with and without augmentation. The distortion here refers to the maximal difference between the unperturbed $\boldsymbol{x}$ and  perturbed $\boldsymbol{x}(\boldsymbol{\xi})$ grids averaged over $1000$ grids used to augment dataset.}
    \label{fig:sensitivity}
    \end{center}
    \vskip -0.2in
\end{figure}

Before training of augmented dataset, it is instructive to evaluate the network trained without augmentation on the augmented train set. This way, we can estimate the degree of equivariance current neural networks have.

More specifically, for this experiment, we take a dataset for stationary-diffusion equation \eqref{eq:stationary diffusion} (Elliptic alpha), and train DilResNet and FNO. After that, we generate a set of augmented datasets using increasingly distorted grids and evaluate neural networks on them. Results are reported in \cref{fig:sensitivity}.

\begin{table*}[!ht]
    \centering
    \begin{tabular}{@{}llcccccc@{}}
        \toprule
         & & \multicolumn{3}{c}{simple datasets} & \multicolumn{3}{c}{complex datasets}\\\cmidrule(r){3-5}\cmidrule(r){6-8}
        Equation & Model & $\times$ & $\surd$ & $g$ & $\times$ & $\surd$ & $g$ \\\cmidrule(r){1-8}
    	   \multirow{6}{*}{Convection-diffusion} &FNO & $0.067$ & $0.048$ & $28\%$ & $0.510$ & $0.418$ & $18\%$ \\
            & DeepONet  &  $0.675$ & $0.567$ & $16\%$ & --- & --- & \\
            & DilResNet & $0.023$ & $0.010$ & $56\%$ & $0.312$ & $0.225$ & $28\%$\\
            
            & MLP & $0.094$ & $0.050$ & $49\%$ & $0.566$ & $0.496$ & $12\%$\\
            & U-Net & $0.069$ & $0.031$ & $55\%$ & $0.419$ & $0.364$ & $13\%$\\
            & SNO & $0.086$ & $0.066$ & $23\%$ & $0.416$ & $0.373$ & $10\%$\\ 
         \cmidrule(r){2-8}
    	\multirow{6}{*}{Elliptic alpha} &FNO & $0.066$ & $0.036$ & $46\%$ & $0.306$ & $0.207$ & $32\%$\\  
        
            & DeepONet  & ---  & $0.826$ & & --- & --- &\\
            
            
            
            & DilResNet  & $0.105$ & $0.021$ & $80\%$ & $0.160$ & $0.133$ & $17\%$\\
            
            & MLP & $0.088$ & $0.053$ & $40\%$ & $0.322$ & $0.253$ & $21\%$\\
            & U-Net & $0.093$ & $0.070$ & $25\%$ & $0.386$ & $0.194$ & $50\%$\\
            & SNO & $0.082$ & $0.050$ & $39\%$ & $0.251$ & $0.209$ & $17\%$\\ 
        \cmidrule(r){2-8}
            \multirow{6}{*}{Elliptic beta} &FNO & $0.034$ & $0.021$ & $38\%$ & $0.181$ & $0.126$ & $30\%$\\ 
         
            & DeepONet  & ---  & $0.832$ & & --- & $0.946$ & \\
            
            
            
            & DilResNet  & $0.099$ & $0.022$ & $78\%$ & $0.089$ & $0.062$ & $30\%$\\
            
            & MLP & $0.069$ & $0.035$ & $50\%$ & $0.238$ & $0.138$ & $42\%$\\
            & U-Net & $0.070$ & $0.067$ & $4\%$ & $0.170$ & $0.143$ & $16\%$\\
            & SNO & $0.068$ & $0.038$ & $44\%$ & $0.187$ & $0.144$ & $23\%$\\ 
        \cmidrule(r){2-8}
            \multirow{6}{*}{Wave} & FNO & $0.200$ & $0.159$ & $21\%$ & $0.650$ & $0.628$ & $3\%$\\  
         
            & DeepONet  & ---  & --- & & --- & --- &\\
            
            
            
            & DilResNet  & $0.053$ & $0.048$ & $9\%$ & $0.43$ & $0.38$ & $12\%$\\
            
            & MLP & $0.313$ & $0.295$ & $6\%$ & --- & $0.99$ &\\
            & U-Net & --- & --- & & $0.57$ & $0.52$ & $9\%$\\
            & SNO & $0.37$ & $0.37$ & $0\%$ & --- & --- & \\\cmidrule(r){2-8}
            & & \multicolumn{3}{c}{$v^{1}$} & \multicolumn{3}{c}{$v^{2}$}\\\cmidrule(r){3-5}\cmidrule(r){6-8}
            \multirow{5}{*}{Navier-Stokes} & FNO & $0.005$ & $0.003$ & $40\%$ & $0.022$ & $0.010$ & $55\%$\\
            & UNet & $0.019$ & $0.09$ & $53\%$ & $0.069$ & $0.037$ & $46\%$\\
            & DilResNet & $0.021$ & $0.015$ & $29\%$ & $0.073$ & $0.045$ & $38\%$\\
            & MLP & $0.082$ & $0.066$ & $38\%$ & $0.082$ & $0.066$ & $20\%$\\
            & SNO & $0.004$ & $0.003$ & $25\%$ &  $0.013$ & $0.008$ & $38\%$\\
        \bottomrule 
    \end{tabular}
    \caption{Relative test errors and gain for $D=2$ datasets. Symbols $\surd$ and $\times$ mark results with and without augmentation respectively. We put --- when network fails to reach test error below $1.0$.}
    \label{table:results_2D}
\end{table*}

As we can see, if distortion is comparable with $10^{-2}$, which is a spacing of the original grid, the neural network can handle the modified dataset quite well. However, the further distortion increase leads to substantial performance deterioration. So for distortion of about $10$ original grid spacing, the trained network is unusable. Networks trained with augmentation retain good relative error even on distorted grids.

Interestingly, FNO can handle distortion slightly better --- four-fold increase against seven-fold for DilResNet. Qualitatively, FNO starts with $4\%$ and ends with $18\%$, whereas DilResNet starts with $10\%$ and ends with $70\%$

We stress that for equations other than elliptic the equivariance is not observed (see results in \cref{table:sensitivity}). Besides, for the elliptic equation, equivariance is confirmed only for the specific distortions of the grid. It is not obvious that the same result holds for other grid transformations, so we do not claim that we achieved general covariance.

\subsection{Statistical study for $D=1$ problems}
Given that neural networks fail to produce correct predictions for augmented dataset, the next natural step is to introduce augmented samples on the training stage. Here we describe relevant experiments for $D=1$ datasets.

In this section we use $m$ to denote augmentation factor. If original train set consists of $N_{\text{train}}$ points, after augmentation with factor $m$ the modified train set has $(1+m)N_{\text{train}}$ points.

For each equation we consider the following parameters: $N_{\text{train}} = 500, 1000, 1500, 2000$; $N_{\text{test}} = 1000$; augmentation factor $m = 1, 2, 3, 4$. For each set of parameters we perform five runs with different seeds controlling network initialization and random grids generated for the augmentation. In \cref{appendix:Supplementray results} one can find relative test errors averaged with respect to these five runs. Here we present only aggregated results.

First, \cref{table:equation vs network} contain results averaged with respect to sizes of train set and augmentation factors. We can see that gain is positive on average for all networks and equations. Note, however, that we observe negative gain, i.e., the augmentation fail to improve test error. This occurres mainly for weak models such as DeepONet and MLP that are othen fail to achieve reasonable test error.

Second, \cref{table:augmentation-test} contain gains averaged over equations and networks. Generally, we observe that augmentation is more helpful for larger datasets. Our best current explanation is that our augmentation procedure is poorly-calibrated, i.e., the augmented data is completely out of distribution. If this is the case, the increase of the train set may improve the overlap between the train data distribution and augmented data distribution.

\subsection{Augmentation for $D=2$ problems}
We report results for a single run for two-dimensional problems in \cref{table:results_2D}.

One can see that, augmentation reduces relative test error for all cases when the network can generalize. It also helps DeepONet to reach test errors smaller than one for the Elliptic alpha and Elliptic beta datasets.

The most intriguing part are the results for the stationary diffusion equation. As we explained before, in $D=2$ we consider two datasets for the elliptic equation: Elliptic alpha and Elliptic beta. Elliptic alpha has a complete set of distinct diffusion coefficients that form a symmetric positive definite matrix.

Contrary to that, Elliptic beta has a single positive diffusion coefficient, so the diffusion matrix $a^{ij}$ in \eqref{eq:stationary diffusion} is proportional to the diagonal matrix. On the other hand, after the augmentation, i.e., for the transformed equation, the Elliptic beta train set has nonzero off-diagonal contributions to the diffusion matrix (see \cref{table:PDEs transforms}). At the same time, we still have a diffusion matrix proportional to the diagonal for the test set. Despite this discrepancy, augmentation still improves the test error.

This result suggests that it can be beneficial to embed a given family of equations into a larger parametric family and perform augmentation for that extended set.

For the Navier-Stokes equation, we can also see that augmentation improves test error for both components of speed $v^{1},\,v^{2}$. This provides evidence that the method is applicable without difficulties to complex geometries. We also tried to train DeepONet, but failed to obtain relative errors comparable to other networks, additional results are available in \cref{appendix:Supplementray results}.

\section{Related research}
We know of two articles directly related to the augmentation techniques for neural PDE solvers \cite{brandstetter2022lie}, \cite{li2022physics}.

The article \cite{li2022physics} is secondary since it is only a minor extension of \cite{brandstetter2022lie}. We do not discuss it further.

In \cite{brandstetter2022lie}, authors consider Lie point symmetries. Namely, to perform augmentation, they use smooth transformations that preserve a solution set of PDE (map a given solution to the other one) and form a group with a structure of the continuous manifold (Lie group). As the authors explain, Lie point symmetries, in a certain sense, provide an exhaustive set of possible transformations suitable for augmentation. Given that, it is appropriate to highlight what distinguishes our research from \cite{brandstetter2022lie}.

In \cite{brandstetter2022lie} symmetries of a \textit{fixed PDE} are used. In place of that, we consider mappings that leave us within \textit{a particular family of PDEs}. Such transformations are more abundant, easier to find, and more suitable for physical systems with spatiotemporal dependencies, e.g., Maxwell equations in macroscopic media \cite{MR0436782}, wave propagation in non-homogeneous media \cite{MR693455}, and fluid flow through porous media \cite{alt1985nonsteady}, e.t.c. In particular, the local distortion of coordinates is a safe choice for a large set of PDEs because of the general covariance principle \cite{post1997formal}.

Among other contributions indirectly related to our approach, we can mention articles on deep learning for PDE that deal with complex geometries \cite{gao2021phygeonet}, \cite{li2022fourier}. The end of this line of research is to generalize physics-informed neural networks and neural operators from rectangular to more general domains. These works also contain particular equations in a transformed form, but for an entirely different reason.

Similar coordinate transformations are abundant in classical scientific computing \cite{MR3675718}, \cite{MR1300634}. On the one hand, the grid defines the geometry of the domain. On the other hand, the grid is refined (h-refinement \cite{li1997adaptive}, \cite{baker1997mesh}) or transported (r-refinement \cite{baker1997mesh}) to improve accuracy (e.g., equidistribution principle \cite{chen1994error}). One can implement both approaches using partial differential equations (elliptic and hyperbolic equations) or analytic mappings (algebraic methods \cite{smith1982algebraic}, \cite{gordon1973construction}). In the present research, we gravitate toward the latter because it is computationally cheap, and derivatives are readily available.

From the broader perspective, the entire subfield of geometric machine learning \cite{bronstein2021geometric} deals with related issues. Typically, the quest is to design a neural network for which the invariance or equivariance holds for a chosen set of transformations. The most relevant works of this sort are \cite{cheng2019covariance}, \cite{weiler2021coordinate}, \cite{wang2020incorporating}. In the first two articles, the authors develop gauge-invariant convolutions (general covariance), but PDEs are not in question, and the generalization to neural operators is not yet available. In the third article, the authors design neural networks that respect selected symmetries of the Navier-Stokes equation.

We want to point out that invariance and equivariance principles are often of no use for PDE problems. First, transformations of the physical fields can fail to be covariant or contravariant, as shown by the example of convection-diffusion equation \cref{table:PDEs transforms}. Second, the symmetries are often not apparent when PDE in question has spatial dependence or is defined in complex geometry. For example, lid-driven cavity flow \cref{eq:cavity flow} breaks all symmetries of the Navier-Stokes equation listed in \cite{wang2020incorporating}.

\section{Conclusion and further research}
We demonstrated how to construct augmentation based on general covariance. The essence of the approach is the observation that it is possible to use change of coordinates to produce novel solutions from the old ones. This is possible because physical phenomena do not depend on the choice of coordinates, so in the new coordinate system, the type of the equation persists, but parameters change. These new parameters along with the solution in new coordinate system can be used as additional train samples.

The proposed augmentation systematically improves test error for all considered architectures. Besides that, it is architecture-agnostic and generalizes well on other equations, especially defined in complex geometries suitable for body-fitted meshes.

A lot of improvements to the proposed approach are possible. The list below sums up a few possibilities:
\begin{enumerate}
    \item \textit{More complex structured grids.} As we showed in the Navier-Stokes example, it is straightforward to extend our approach to situations where there is an intermediate mapping from the physical domain to the computational domain. The considered example is elementary, so it is desirable to test augmentation on more challenging problems. It is also interesting to consider our augmentation approach with the architectures that already contain mapping as part of the network, e.g., \cite{gao2021phygeonet}, \cite{li2022fourier}.
        
    \item \textit{Adaptive augmentation.} In the present research, we generate random grids to perform augmentation. It should be more advantageous to actively select grids based on the neural network performance.
    
    \item \textit{Unstructured grids.} Our approach, as described here, is, in principle, applicable to the unstructured grids. The main problem is that it is not obvious how to construct a mapping and generate a deformed grid such that it is still acceptable from the computational perspective.
    
    \item \textit{Covariant neural operators.} It would be interesting to adapt or generalize results from \cite{weiler2021coordinate}, \cite{wang2020incorporating} to construct covariant neural operators. Less ambitious, it is possible to improve the training protocols for neural networks using the Jacobi matrix and determinant as input features. This way, one escapes the need to introduce higher derivatives in transformed equations \cref{table:PDEs transforms}.
        
    \item \textit{Transformations between parametric families.} A very interesting feature that we observe is that augmentation still helps even when it is performed for a larger parametric family of equations than needed. One possible extension is to dispense with coordinate covariance and consider more general transformations that map one family of PDEs to another family.
    
    \item \textit{Time-dependent coordinate transformations.} We only consider spatial transformations. It is possible to use time-dependent transformations. For example, one can construct a grid, that is deformed at $t=0$ and approaches its non-deformed state when $t$ increases.
\end{enumerate}
\section{Acknowledgements}
The work was supported by the Analytical center under the RF Government (subsidy agreement 000000D730321P5Q0002, Grant No. 70-2021-00145 02.11.2021)
\bibliography{refs}
\bibliographystyle{icml2023}
\newpage
\appendix
\onecolumn
\section{Coordinate transformations}
\label{appendix:Coordinate transformations}
In this appendix we collect standard material related to coordinate transformations and transformation laws for differential operators.

The material of this section is available in many sources \cite{MR3675718}, \cite{MR1441306}, \cite{MR1247707} and presented here merely for convenience of the reader.

Everywhere in this section the Einstein's summation notation is used, e.g., $a_{\alpha}b^{\alpha} \equiv \sum_{\alpha} a_{\alpha}b^{\alpha}$, and as coordinate transformations we consider \eqref{eq:coordinate transformation}.
\subsection{Some relations for the first and the second derivatives}
\label{appendix:relations}
In this section we explain several identities we find useful for expressing equations in the covariant form. The main problem is that when we have a mapping $\boldsymbol{x}(\boldsymbol{\xi})$ (defined below) given in some explicit form and have no closed-form expression for the inverse mapping, it is inconvenient to use derivatives with respect to $\boldsymbol{x}$. Since such derivatives appear in plenitude in PDEs after coordinate transformation, our chief goal is to find relations that allows us to rewrite them using derivatives with respect to $\boldsymbol{\xi}$.

We use Jacobi's formula for differentiable matrix-valued function $\boldsymbol{A}(t)$ without a proof:
\begin{equation}
    \label{eq:Jacobi_formula}
    \frac{d}{dt}\det \boldsymbol{A}(t) = \det \boldsymbol{A}(t)\text{ tr} \left(\boldsymbol{A}^{-1}(t)\frac{d \boldsymbol{A}(t)}{dt}\right).
\end{equation}
The short note \cite{MR315183} contains a concise derivation.

We start with the relation between first derivatives
\begin{equation}
    \label{eq:matrix_vs_inverse}
    \frac{\partial x^{i}}{\partial \xi^{\alpha}} \frac{\partial \xi^{\alpha}}{\partial x^{j}} = \delta_{ij}
    =
    \begin{cases}
        1,\,i = j;\\
        0,\,i \neq j,
    \end{cases}
\end{equation}
which follows from chain rule applied to the function $\boldsymbol{x}(\boldsymbol{\xi}(\boldsymbol{x})) = \boldsymbol{x}$:
\begin{equation}
    \delta_{ij} = \frac{\partial x^{i}}{\partial x^{j}} = \frac{\partial x^{i}(\xi^{1}(x^{1}, \dots, x^{D}),\dots,\xi^{D}(x^{1}, \dots, x^{D}))}{\partial x^{j}} = \frac{\partial x^{i}}{\partial \xi^{\alpha}}\frac{\partial \xi^{\alpha}}{\partial x^{j}}.
\end{equation}
It is also convenient to rewrite identity \cref{eq:matrix_vs_inverse} in a matrix form
\begin{equation}
    \label{eq:matrix_vs_inverse_m}
    \mathcal{J}_{i\alpha} \equiv \frac{\partial x^{i}}{\partial \xi^{\alpha}},\, \boldsymbol{\mathcal{J}}\boldsymbol{\mathcal{J}}^{-1} = \boldsymbol{I},
\end{equation}
where $\boldsymbol{\mathcal{J}}$ is Jacobi matrix and $\boldsymbol{I}$ is the identity matrix.

Next identity we consider is the equation for the derivative of $J=\det \boldsymbol{\mathcal{J}}$, i.e., the determinant of Jacobi matrix:
\begin{equation}
    \label{eq:d_Jacobi_det}
    \frac{\partial}{\partial \xi^{k}} J = J\frac{\partial \xi^{m}}{\partial x^{i}} \frac{\partial^2 x^{i}}{\partial \xi^{m}\partial \xi^{k}}.
\end{equation}
The equation immediately follows from Jacobi's equation \eqref{eq:Jacobi_formula} and the definition of inverse \eqref{eq:matrix_vs_inverse}, \eqref{eq:matrix_vs_inverse_m} for Jacobi matrix.

Another relation useful in derivations reads
\begin{equation}
    \label{eq:second_derivative}
    \frac{\partial^2 \xi^{\alpha}}{\partial x^{i}\partial x^{j}} = -\frac{\partial \xi^{\gamma}}{\partial x^{i}}\frac{\partial \xi^{\theta}}{\partial x^{j}} \frac{\partial \xi^{\alpha}}{\partial x^{k}} \frac{\partial^2 x^{k}}{\partial \xi^{\gamma}\partial \xi^{\theta}}.
\end{equation}
To prove the relation, we find a derivative of \eqref{eq:matrix_vs_inverse} as follows
\begin{equation}
    \label{eq:second_derivative_d}
    0 = \frac{\partial}{\partial x^{k}}\left(\frac{\partial x^{i}}{\partial \xi^{\alpha}} \frac{\partial \xi^{\alpha}}{\partial x^{j}}\right) = \frac{\partial}{\partial x^{k}}\left(\frac{\partial x^{i}}{\partial \xi^{\alpha}}\right) \frac{\partial \xi^{\alpha}}{\partial x^{j}} + \frac{\partial x^{i}}{\partial \xi^{\alpha}}\frac{\partial^2 \xi^{\alpha}}{\partial x^{j} \partial x^{k}} = \frac{\partial \xi^{\beta}}{\partial x^{k}}\frac{\partial^2 x^{i}}{\partial \xi^{\alpha} \partial \xi^{\beta}} \frac{\partial \xi^{\alpha}}{\partial x^{j}} + \frac{\partial x^{i}}{\partial \xi^{\alpha}}\frac{\partial^2 \xi^{\alpha}}{\partial x^{j} \partial x^{k}},
\end{equation}
and multiply by the inverse Jacobi matrix.

The last identity we need reads
\begin{equation}
    \label{eq:div_0}
    \frac{1}{J}\frac{\partial}{\partial \xi^{j}} \left(J \frac{\partial \xi^{j}}{\partial x^{i}}\right) = 0.
\end{equation}
To prove \eqref{eq:div_0} we apply \eqref{eq:d_Jacobi_det} (Jacobi formula) and obtain
\begin{equation}
    \frac{1}{J}\frac{\partial}{\partial \xi^{\alpha}} \left(J \frac{\partial \xi^{\alpha}}{\partial x^{j}}\right) = \frac{\partial \xi^{\beta}}{\partial x^{k}}\frac{\partial^2 x^{k}}{\partial \xi^{\alpha} \partial \xi^{\beta}} \frac{\partial \xi^{\alpha}}{\partial x^{j}} + \frac{\partial}{\partial \xi^{\alpha}}\frac{\partial \xi^{\alpha}}{\partial x^{j}} = \frac{\partial \xi^{\beta}}{\partial x^{k}}\frac{\partial^2 x^{k}}{\partial \xi^{\alpha} \partial \xi^{\beta}} \frac{\partial \xi^{\alpha}}{\partial x^{j}} + \frac{\partial x^{k}}{\partial \xi^{\alpha}}\frac{\partial^2 \xi^{\alpha}}{\partial x^{j} \partial x^{k}}.
\end{equation}
The last expression is zero since it is a particular form of more general identity \eqref{eq:second_derivative_d} with $i=k$.
\subsection{Selected differential operators under coordinate transformation}
\label{appendix:differential operators}
Using results from \cref{appendix:relations} we derive transformation laws for particular differential operators.

Two equations for the first derivative
\begin{equation}
    \label{eq:first_non_conservative}
    c^{j}\frac{\partial\phi}{\partial x^{j}} = c^{j}\frac{\partial \xi^{\alpha}}{\partial x^{j}}\frac{\partial\phi}{\partial \xi^{\alpha}}
\end{equation}
and for the second
\begin{equation}
    \label{eq:second_non_conservative}
    a^{kj} \frac{\partial^2\phi}{\partial x^{j}\partial x^{k}} = a^{kj}\frac{\partial \xi^{\beta}}{\partial x^{j}} \frac{\partial}{\partial \xi^{\beta}}\left(\frac{\partial \xi^{\gamma}}{\partial x^{k}}\frac{\partial \phi}{\partial \xi^{\gamma}}\right) = a^{kj} \frac{\partial \xi^{\beta}}{\partial x^{j}} \frac{\partial \xi^{\gamma}}{\partial x^{k}}\frac{\partial^2 \phi}{\partial \xi^{\beta}\partial \xi ^{\gamma}} + a^{kj}\frac{\partial^2 \xi^{\gamma}}{\partial x^{k}\partial x^{j}} \frac{\partial \phi}{\partial \xi^{\gamma}}
\end{equation}
follow simply from the chain rule. To apply these equations when only a mapping $\boldsymbol{x}(\boldsymbol{\xi})$ is known one needs to use relations \eqref{eq:matrix_vs_inverse} and \eqref{eq:second_derivative} for \eqref{eq:second_non_conservative}.

Conservative forms of the equations above read
\begin{equation}
    \label{eq:first_conservative}
    \frac{\partial}{\partial x^{\alpha}}\left(c^{\alpha} \phi\right) = \frac{1}{J}\frac{\partial}{\partial \xi^{k}}\left(Jc^{\alpha}\frac{\partial \xi^{k}}{\partial x^{\alpha}}\phi\right)
\end{equation}
and
\begin{equation}
    \label{eq:second_conservative}
    \frac{\partial}{\partial x^{k}}\left(a^{kj}\frac{\partial\phi}{\partial x^{j}}\right) = \frac{1}{J}\frac{\partial}{\partial \xi^{k}}\left(J\left( a^{\alpha j}\frac{\partial \xi ^{k}}{\partial x^{\alpha}}  \frac{\partial \xi^{\beta}}{\partial x^{j}} \right)\frac{\partial\phi}{\partial \xi^{\beta}}\right)
\end{equation}
respectively.

\cref{eq:first_conservative} is straightforward to confirm using \eqref{eq:div_0}. Indeed,
\begin{equation}
    \frac{1}{J}\frac{\partial}{\partial \xi^{k}}\left(Jc^{\alpha}\frac{\partial \xi^{k}}{\partial x^{\alpha}}\phi\right) = \underbrace{\frac{1}{J}\frac{\partial}{\partial \xi^{k}}\left(J\frac{\partial \xi^{k}}{\partial x^{\alpha}}\right)c^{\alpha}\phi}_{=0} + \frac{\partial \xi^{k}}{\partial x^{\alpha}}\phi\frac{\partial c^{\alpha}}{\partial \xi^{k}} + \frac{\partial \xi^{k}}{\partial x^{\alpha}}c^{\alpha}\frac{\partial \phi}{\partial \xi^{k}} = \frac{\partial}{\partial x^{\alpha}}\left(c^{\alpha} \phi\right).
\end{equation}

\cref{eq:second_conservative} is easier to derive from the analogous result for the divergence of a vectors field
\begin{equation}
    \label{eq:divergence}
    \frac{\partial}{\partial x^{k}} f^{k}(x) = \frac{1}{J}\frac{\partial}{\partial \xi^{j}} \left(J f^{i}(x(\xi)) \frac{\partial \xi^{j}}{\partial x^{i}}\right).
\end{equation}
\cref{eq:divergence} itself trivially follows from \eqref{eq:div_0}.

Having \eqref{eq:divergence} we derive \eqref{eq:second_conservative} as follows
\begin{equation}
    \frac{\partial}{\partial x^{k}}\left(a^{kj}\frac{\partial\phi}{\partial x^{j}}\right) = \frac{1}{J}\frac{\partial}{\partial \xi^{k}} \left(J a^{\alpha j}\frac{\partial\phi}{\partial x^{j}} \frac{\partial \xi^{k}}{\partial x^{\alpha}}\right) = \frac{1}{J}\frac{\partial}{\partial \xi^{k}}\left(J\left( a^{\alpha j}\frac{\partial \xi ^{k}}{\partial x^{\alpha}}  \frac{\partial \xi^{\beta}}{\partial x^{j}} \right)\frac{\partial\phi}{\partial \xi^{\beta}}\right).
\end{equation}
\subsection{Selected PDEs under coordinate transformation}
\label{appendix:equations}
Results from \cref{appendix:differential operators} allows to derive transformed form of a large set of equations. We confine our attention to stationary diffusion \eqref{eq:stationary diffusion}, convection-diffusion \eqref{eq:convection diffusion} and wave equations \eqref{eq:wave equation}.

We start with the derivation of the transformation law for \eqref{eq:stationary diffusion}. According to \eqref{eq:second_conservative}, after the transformation \eqref{eq:coordinate transformation} equation becomes
\begin{equation}
    \label{eq:stationary diffusion transformed}
    \frac{\partial}{\partial \xi^{k}}\left(J\left( a^{\alpha j}(\boldsymbol{x}(\boldsymbol{\xi}))\frac{\partial \xi ^{k}}{\partial x^{\alpha}}  \frac{\partial \xi^{\beta}}{\partial x^{j}} \right)\frac{\partial u(\boldsymbol{x}(\boldsymbol{\xi}))}{\partial \xi^{\beta}}\right) = Jf(\boldsymbol{x}(\boldsymbol{\xi})),
\end{equation}
with the same boundary conditions. So, it is exactly the same parametric equation as \eqref{eq:stationary diffusion} but with different parameters.

Next, we consider the convection-diffusion equation \eqref{eq:convection diffusion}. Using again \eqref{eq:second_conservative} and \eqref{eq:first_conservative} we obtain a transformed form
\begin{equation}
    \label{eq:convection diffusion transformed 1}
    \frac{\partial}{\partial t}\phi(\boldsymbol{x}(\boldsymbol{\xi}), t) + \frac{1}{J} \frac{\partial}{\partial \xi^{i}} \left(Jv^{k}(\boldsymbol{x}(\boldsymbol{\xi}))\frac{\partial \xi^{i}}{\partial x^{k}}\phi(\boldsymbol{x}(\boldsymbol{\xi}), t)\right) = \frac{1}{J}\frac{\partial}{\partial \xi^{k}}\left(J\left( a^{\alpha j}(\boldsymbol{x}(\boldsymbol{\xi}))\frac{\partial \xi ^{k}}{\partial x^{\alpha}}  \frac{\partial \xi^{\beta}}{\partial x^{j}} \right)\frac{\partial\phi(\boldsymbol{x}(\boldsymbol{\xi}), t)}{\partial \xi^{\beta}}\right).
\end{equation}

This time we can see that \eqref{eq:convection diffusion transformed 1} does not have the same parametric form as \eqref{eq:convection diffusion}. To resolve the problem we define a new field $\psi(\boldsymbol{\xi}, t) \equiv J \phi(\boldsymbol{x}(\boldsymbol{\xi}), t)$. Multiplying both sides of \cref{eq:convection diffusion transformed 1} by $J$ we observe that left hand side has a desired parametric form. For the right hand side we proceed as follows
\begin{equation}
    \begin{split}
    \frac{\partial}{\partial \xi^{k}}\left(Ja^{\alpha j}(\boldsymbol{x}(\boldsymbol{\xi}))\frac{\partial \xi ^{k}}{\partial x^{\alpha}}  \frac{\partial \xi^{\beta}}{\partial x^{j}} \frac{\partial}{\partial \xi^{\beta}}\left(\frac{\psi(\boldsymbol{x}(\boldsymbol{\xi}), t)}{J}\right)\right) &= \frac{\partial}{\partial \xi^{k}}\left(a^{\alpha j}(\boldsymbol{x}(\boldsymbol{\xi}))\frac{\partial \xi ^{k}}{\partial x^{\alpha}}  \frac{\partial \xi^{\beta}}{\partial x^{j}} \frac{\partial\psi(\boldsymbol{x}(\boldsymbol{\xi}), t)}{\partial \xi^{\beta}}\right)\\
    &- \frac{\partial}{\partial \xi^{k}}\left(a^{\alpha j}(\boldsymbol{x}(\boldsymbol{\xi}))\frac{\partial \xi ^{k}}{\partial x^{\alpha}}  \frac{\partial \xi^{\beta}}{\partial x^{j}} \frac{\partial \xi^{\gamma}}{\partial x^{\rho}}\frac{\partial^2 x^{\rho}}{\partial \xi^{\gamma}\partial \xi^{\beta}}\psi(\boldsymbol{x}(\boldsymbol{\xi}))\right),
    \end{split}
\end{equation}
where we used \eqref{eq:d_Jacobi_det}. As we can see the right hand side of \eqref{eq:convection diffusion transformed 1} introduces additional contribution to the convection term, and with this contribution the final equation has the same parametric form as the original convection-diffusion equation \eqref{eq:convection diffusion}.

In the case of two-way wave equation \eqref{eq:wave equation} we apply \eqref{eq:first_non_conservative} and \eqref{eq:second_non_conservative} to obtain
\begin{equation}
    \label{eq:wave equation transformed}
    \begin{split}
        \frac{\partial^2 \rho(\boldsymbol{x}(\boldsymbol{\xi}), t)}{\partial t^2} + v^{i}(\boldsymbol{x}(\boldsymbol{\xi}))\frac{\partial \xi^{\alpha}}{\partial x^{i}} \frac{\partial \rho(\boldsymbol{x}(\boldsymbol{\xi}), t)}{\partial \xi^{\alpha}} = c^{kj}(\boldsymbol{x}(\boldsymbol{\xi})) \frac{\partial \xi^{\beta}}{\partial x^{j}} \frac{\partial \xi^{\gamma}}{\partial x^{k}}\frac{\partial^2 \rho(\boldsymbol{x}(\boldsymbol{\xi}), t)}{\partial \xi^{\beta}\partial \xi ^{\gamma}} &+ c^{kj}(\boldsymbol{x}(\boldsymbol{\xi}))\frac{\partial^2 \xi^{\gamma}}{\partial x^{k}\partial x^{j}} \frac{\partial \rho(\boldsymbol{x}(\boldsymbol{\xi}), t)}{\partial \xi^{\gamma}} \\&+ e(\boldsymbol{x}(\boldsymbol{\xi})\rho(\boldsymbol{x}(\boldsymbol{\xi}).
    \end{split}
\end{equation}
All transformation laws we derived can be found in \cref{table:PDEs transforms}.
\section{Coordinate derivatives for transfinite interpolation}
\label{appendix:transfinite interpolation}
For the case of transfinite interpolation \cref{eq:2D mapping}, Jacobi matrix and higher order derivatives slightly simplify. Here we provide explicit expressions for the $D=2$ case:

\begin{enumerate}
    \item Jacobi matrix
    
    \begin{equation}
        \mathcal{J}_{ij} = \frac{\partial x^{i}}{\partial \xi^{j}} = 
        \begin{pmatrix}
            y_1^{'}(\xi^{1})(1-\xi^{2}) + y_2^{'}(\xi^{1})\xi^{2} & y_2(\xi^{1}) - y_1(\xi^{1})\\
            y_4(\xi^{2}) - y_3(\xi^{2}) & (1-\xi^{1})y_3^{'}(\xi^{2}) + \xi^{1} y_4^{'}(\xi^{2})\\
        \end{pmatrix}.
    \end{equation}

    \item Inverse Jacobi matrix

        \begin{equation}
            \left(\mathcal{J}^{-1}\right)_{ij} = \frac{\partial \xi^{i}}{\partial x^{j}} = 
            \frac{1}{J}\begin{pmatrix}
                (1-\xi^{1})y_3^{'}(\xi^{2}) + \xi^{1} y_4^{'}(\xi^{2}) & y_1(\xi^{1}) - y_2(\xi^{1})\\
                y_3(\xi^{2}) - y_4(\xi^{2}) & y_1^{'}(\xi^{1})(1-\xi^{2}) + y_2^{'}(\xi^{1})\xi^{2}\\
            \end{pmatrix},\,J = \det \mathcal{J}.
        \end{equation}

    \item Second derivatives
    
        \begin{equation}
            \begin{split}
            &\frac{\partial^2 x^{1}}{\partial \xi^{i} \partial \xi^{j}} = 
            \begin{pmatrix}
                y_1^{''}(\xi^{1})(1-\xi^{2}) + y_2^{''}(\xi^{1})\xi^{2} & P^{'}_2(\xi^{1}) - P^{'}_1(\xi^{1})\\
                P^{'}_2(\xi^{1}) - P^{'}_1(\xi^{1}) & 0\\
            \end{pmatrix},\\
            &\frac{\partial^2 x^{2}}{\partial \xi^{i} \partial \xi^{j}} = 
            \begin{pmatrix}
                0 & P^{'}_4(\xi^{2}) - P^{'}_3(\xi^{2})\\
                P^{'}_4(\xi^{2}) - P^{'}_3(\xi^{2}) & (1-\xi^{1})y_3^{''}(\xi^{2}) + \xi^{1} y_4^{''}(\xi^{2})\\
            \end{pmatrix}.
            \end{split}
        \end{equation}
\end{enumerate}
\section{Architectures and training details}
\label{appendix:Architectures and training details}
In this section, we provide extended comments on the architectures used and collect in \cref{table:training_details} the description of the optimization process.

\begin{table}[b]
    \centering
    \caption{Training details: $\nu$ --- learning rate, $\nu$ decay / epoch --- weight decay per epoch, $N_{\text{epoch}}$ --- number of epoch used for training, $N_{\text{batch}}$ --- batch size, $N_{\text{params}}$ --- number of network parameters.
    \\ ${}^{\dagger}$ in DeepONet, inverse time function means $\nu=\nu / (1 + 0.5 \cdot \text{steps} / g)$, where $g=\lfloor \text{epoch}/5 \rfloor$
    \\ ${}^{\ddagger}$ $1000$ epochs were used for $D=2$ datasets.}
    \vskip 0.15in
    \begin{tabular}{@{}lccccccc}
        \toprule
        Network & $\nu$ & $\nu$ decay / epoch & weight decay & $N_{\text{epoch}}$ & $N_{\text{batch}}$ & \multicolumn{2}{c}{$N_{\text{params}}$} \\ \cmidrule{2-8}
        & & & & & & $D=1$ & $D=2$\\ \cmidrule{7-8}
        FNO & $10^{-3}$ & --- & $10^{-4}$  &  $500$ &  $200$ & $549\times10^3$ & $236 \times 10^{4}$\\
        DeepONet & $10^{-3}$ & --- &  inverse time$^{\dagger}$ & $2 \times 10^{4}$ & full train set &  $150\times10^3$  & $816 \times 10^{4}$\\
        rFNO & $10^{-3}$ & $0.5 \big/ 100$ & $10^{-2}$ & $500$ & $30$ & $215 \times 10^{3}$ & --- \\
        MLP & $10^{-3}$ & $0.5 \big/ 100$ & $10^{-2}$ & $500^{\ddagger}$ & $30$ & $63 \times 10^{3}$ & $113\times 10^3$\\
        DilResNet & $10^{-3}$ & $0.5 \big/ 100$ & $10^{-2}$ & $500$ & $30$ & $87 \times 10^{3}$ & $261\times 10^{3}$\\
        U-Net & $10^{-3}$ & $0.5 \big/ 100$ & $10^{-2}$ & $500$ & $30$ & --- & $263 \times 10^{3}$\\
        SNO  & $10^{-3}$ & $0.5 \big/ 100$ & $10^{-2}$ & $500$ & $30$ & --- & $115 \times 10^{3}$\\
    \bottomrule
    \end{tabular}
    \label{table:training_details}
    \vskip -0.1in
\end{table}

\subsection{rFNO}
rFNO is a variant of FNO \cite{li2020fourier} with two differences.

First, FFT is replaced with projection on the set of trigonometric functions
\begin{equation}
    \label{eq:trigonometric space}
    \mathcal{T}_{N} \equiv \left\{1, \cos\left(2\pi  x\right), \sin\left(2\pi x\right), \cos\left(2\pi  2 x\right), \sin\left(2\pi 2x\right), \dots, \cos\left(2\pi N x\right), \sin\left(2\pi N x\right) \right\}.
\end{equation}
That is, for the input function $f(x)$, the function after transform reads
\begin{equation}
    c_k = \int_{0}^{1} dx\, \phi_{k}(x) f(x),\,\phi_{k}\in\mathcal{T}_{N},
\end{equation}
where the integral is approximated using the trapezoidal rule.

The inverse transform is computed simply as a sum
\begin{equation}
    g(x) = \sum_{k} c_k \phi_k(x),\,\phi_{k}\in\mathcal{T}_{N}.
\end{equation}

Second, between the transformations, FNO \cite{li2020fourier} uses a diagonal tensor, so the resulting matrix performs convolution (see \cite{rippel2015spectral}). In our case we apply a series of convolutions $N_{\text{conv}}$ without activations.

In all $1D$ experiments, we keep $N=8$ modes. Use $N_{\text{conv}}=4$ convolutions with kernel size $3$ in the Fourier space \eqref{eq:trigonometric space}. The encoder lifts the input to the space with $N_{\text{features}}=64$ features. The number of Fourier layers is $4$. We use $\text{ReLU}$ activation functions.

\subsection{MLP}

Since we process functions sampled on the grid, we use MLP that applies a linear operator to each dimension and feature space separately, i.e., if we have a tensor $t_{ijk}$ as an input, the linear (affine transform) layer transforms it as follows
\begin{equation}
    \widetilde{t}_{abc} = \sigma\left(\sum_{ijk} A_{ai} B_{bj} C_{ck}t_{ijk} + e_{abc}\right),
\end{equation}
where $A$, $B$, $C$ are parameters of linear transform and $e$ is bias and $\sigma$ is the activation function.

We use in total $4$ such layers both in $D=1$ and $D=2$, and with the same number of spatial points as input has, and the number of features increased to $64$. Again, $\text{ReLU}$ activation functions were used.

\subsection{DilResNet}

The dilated residual network follows the publication \cite{stachenfeld2021learned}. We use $4$ blocks and $32$ features, each block consists of convolutions with strides $[1, 2, 4, 8, 4, 2, 1]$ each with kernel size $3$. After each block, we put a skip connection. As before, $\text{ReLU}$ activation functions were used.

\subsection{U-Net}

The usual form of U-Net was used \cite{ronneberger2015u}.

U-Net induces the series of grids (levels) each having roughly $\times2$ fewer points, and the number of features doubles. We used $3$ convolutions on each level before $\max$ pooling, transposed convolution for upsampling, and $3$ convolutions for each level after the upsampling. In total we have $4$ layers and start with $10$ features. Again, $\text{ReLU}$ activation functions were used.

\subsection{DeepONet}
DeepONet \cite{lu2019deeponet} consists of two sub-networks, one for encoding the input function $v$ at a fixed number of sensors $x_i, i = 1, \ldots, m$ (branch net), and the other for encoding the locations $\xi$ for the output functions (trunk net). The output of the network can be expressed as
\begin{equation}
\mathcal{G}(v)(\xi)=\sum_{k=1}^p b_k(v) t_k(\xi)+b_0,
\end{equation}
where $b_0$ is the bias, $\{b_k\}_{k=1}^p$ are the outputs of the branch net, and $\{t_k\}_{k=1}^p$ are the outputs of the trunk net.

For the $D=1$ problem, both branch net and trunk net are fully connected neural networks(FNNs). For branch net, we use $4$ layers with $N_{\text{features}}=128$ features; for trunk net, we use $3$ layers with $N_{\text{features}}=128$ features. Additionally, we utilized $\tanh$ as the activation function, and Glorot normal initializer to initialize the weights of DeepONet.

For the $D=2$ problem, the trunk net is fully connected neural networks(FNNs) with $3$ layers, $N_{\text{features}}=128$ features each. For the branch net, we utilized $2$ convolutional layers with $N_{\text{features}}=64$ features and $N_{\text{features}}=128$ features respectively, then we had a flatten layer and two fully connected layers with $N_{\text{features}}=128$ features. Additionally, we utilized $\text{ReLU}$ as the activation function, and Glorot normal initializer to initialize the weights of DeepONet.

\subsection{FNO}
The original form of FNO was proposed in \cite{li2020fourier}.
This network consists of encoder, several Fourier layers and decoder. The $(l+1)$-th Fourier layer can be expressed as
\begin{equation}
z_{l+1}=\sigma\left(\mathcal{F}^{-1}\left(R_l \cdot \mathcal{F}\left(z_l\right)\right)+ \text{conv}(z_l)_{W_l}+b_l\right),
\end{equation}

where $R_l, W_l$ are the weight matrices, $\sigma$ is the activation function, $b_l$ is the bias, and $\mathcal{F}$ is the Fast Fourier transform, $\mathcal{F}^{-1}$ is the inverse, and $\text{conv}$ stands for convolution with kernel size $1$.

For the $D=1$ problem, the number of Fourier layers is $4$. Each layer had $N_{\text{features}}=64$ features and $N=16$ modes. Additionally, we utilized $\text{GELU}$ as the activation function. In addition, we had a linear encoder and decoder implemented as fully connected layers. The encoder was a single layer that lifts the input function to the space with $64$ features. The decoder had two fully connected layers that change the number of features to $64$, then to $128$, and finally to $1$, which is the target number of features for all datasets used.

For the $D=2$ problem, the same number of FNO layers was used. Each layer had $N_{\text{features}}=32$ features and $N=12$ modes. The encoder lifted the input function to the space with $32$ features. The decoder consisted of two fully connected layers, that firstly increased number of features to $128$ and then reduced them to the target number of features. Also, we used $\text{GELU}$ activation functions.

\subsection{SNO}
SNO follows the design of FNO with two notable differences. First, Discrete Cosine Transform replaces FFT. As explained in \cite{fanaskov2022spectral}, this corresponds to the approximate recovery of coefficient in the Chebyshev series, should the input function is sampled on the Chebyshev grid. Second, after DCT, we use three convolutions with kernel size $3$ as a linear kernel. The number of features used in the processor is $32$, number of modes left after truncation is $16$. In total, network contains $4$ layers with the integral operator, ReLU nonlinearities in-between them, linear encoder, and decoder.
\section{Datasets}
\label{appendix:Datasets}
In this section, we provide more details on the dataset generation process. Datasets can be downloaded from \url{https://disk.yandex.ru/d/ArC6jT3TZcKncw}.

\subsection{$D=1$}
To generate input data for $D=1$ PDEs we use two families of trigonometric functions
\begin{equation}
    \label{eq:random diffusion coefficient}
    f_{N}(x) = \sum_{k=0}^{N} c_{k} \cos(2 \pi k x + p_{k})
\end{equation}
and
\begin{equation}
    \label{eq:random rhs}
    g_{N}(x) = \sum_{k=0}^{N} c_{k} \sin(\pi (k+1) x).
\end{equation}
Functions $g_{N}(x)$ with $c_{k}$ from standard normal distribution were used to generate initial conditions for wave and convection-diffusion equations.

Functions $f_{N}(x)$ were used for two purposes. First, we generated positive diffusion coefficients taking $p_{k}, c_{k}, k>0$ from standard normal distribution and fixing $p_{0}=0, c_{0} = \sum_{k>1}\left|c_{k}\right| + \epsilon$ with $\epsilon = 10^{-2}$. Second, we used them with $p_{k}, c_{k}, k\geq0$ taken from standard normal distribution to generate convection coefficient $v(x)$ for the convection-diffusion problem and right-hand side $f(x)$ for the elliptic problem.

The list of generated datasets is as follows:
\begin{enumerate}
    \item Elliptic \eqref{eq:stationary diffusion}
    
        The right-hand side was generated using \eqref{eq:random rhs} with $N=3$ coefficients, for diffusion coefficient we use \eqref{eq:random diffusion coefficient} with $N=5$. We perform discretization with standard FEM method with hat functions on the uniform grid with $100$ points. 
        
    \item Convection-diffusion \eqref{eq:convection diffusion}
    
        Convection coefficient was generated using \eqref{eq:random diffusion coefficient} with $N=5$ coefficients, for diffusion coefficient we use \eqref{eq:random diffusion coefficient} with $N=5$, initial conditions were generated using \eqref{eq:random rhs} with $N=10$. After generation diffusion and convection coefficients are multiplied by $s=0.01$. Spatial discretization is the same as for the Elliptic dataset, for the time-marching Crank-Nicolson scheme was used, final $t=1.0$, $200$ points were used along $t$.
        
    \item Wave$(5)$ \eqref{eq:wave equation}

        All functions were sampled with $N=5$ modes. For initial conditions, we used \eqref{eq:random rhs}, and for source terms, diffusion, and convection coefficients we used \eqref{eq:random diffusion coefficient} multiplied by $s=0.1$. To make the diffusion coefficient positive we squared \eqref{eq:random diffusion coefficient}. For all functions, we used an additional factor for $1/k^2$ for coefficient $k$ to obtain smooth functions. Along the spacial dimension, we use $100$ points and the standard second-order finite-difference discretization, along temporal $1000$ points. As a marching scheme, we used leapfrog and take $t=1.0$ as a final time.
    
    \item Wave$(10)$ \eqref{eq:wave equation}

        The same as previous but with $N=10$ for all sampled functions.
\end{enumerate}

\subsection{$D=2$}
In this case, we use a single family of random functions
\begin{equation}
    \label{eq:random 2D}
    f(x) = \mathcal{R}\left(\sum_{m=-M}^{M}\sum_{n=-M}^{M} e^{2\pi i(mx + ny)}c_{mn}\right),
\end{equation}
where $\mathcal{R}$ is a real part and $c_{mn} = x_{mn} + iy_{mn}$ where both $x_{mn}$ and $y_{mn}$ are samples from standard normal distribution.

The elliptic part of differential operators requires a uniformly positive definite matrix. To enforce this condition we generate Cholesky decomposition of diffusion coefficient $A$ as follows
\begin{equation}
    \label{eq:Cholesky}
    A = I + L L^T,
\end{equation}
where each non-negative element of upper-triangular matrix $L$ is an independently generated function \eqref{eq:random 2D}.

The list of generated datasets is as follows:
\begin{enumerate}
    \item Elliptic alpha \eqref{eq:stationary diffusion}
        Diffusion coefficients were generated in the form of Cholesky decomposition \eqref{eq:Cholesky}. The right-hand side was generated using \eqref{eq:random 2D}. For the simple dataset coefficients in \eqref{eq:random 2D} were multiplied by a factor $s=0.1$ for the complex dataset by the factor $s=0.5$. For both simple and complex datasets we used $M=5$ in \eqref{eq:random 2D}. Bilinear FEM discretization on the rectangular grid with $100\times 100$ points was used for discretization.
    \item Elliptic beta \eqref{eq:stationary diffusion}
        For this dataset the matrix $A$ is diagonal. We generated everything the same way as for the previous datasets, and afterward, drop non-diagonal elements of $A$ and replace $A_{22}$ with $A_{11}$.
    \item Convection-diffusion \eqref{eq:convection diffusion}
        The diffusion coefficients were generated the same way as for the Elliptic alpha dataset. Initial conditions and convection coefficients are taken from \eqref{eq:random 2D}. For each function, $M=5$ is used. For the complex dataset, we rescale coefficients by a factor $s=0.5$ and take $t=1e-2$. For the simple dataset, these parameters are $s=0.1$, $t=1e-2$. As in the $D=1$ case Crank-Nicolson time-marching scheme is used. Spatial discretization is the same as for the Elliptic alpha dataset, $100$ points along a temporal dimension are used.
    \item Wave \eqref{eq:wave equation}
        All random functions were generated the same way as for the convection-diffusion equation. The source term is not used for $D=2$. For the complex dataset, we have $s=0.2$, $t=1$, and for the simple, we put $s=0.2$, $t=1e-1$.
\end{enumerate}
\section{Details on the solution method for lid-driven cavity flow}
\label{appendix:Details on lid-driven cavity flow}
To integrate the equation we use the Chorin projection method that works in three steps:
\begin{enumerate}
    \item Advance speed neglecting pressure term
        \begin{equation}
          \frac{u^{i} - v^{i}_{(n)}}{\Delta t} = \frac{\partial}{\partial x^{k}}\left(-v^{k}_{(n)} v^{i}_{(n)} + \nu \frac{\partial v^{i}_{(n)}}{\partial x^k}\right)
        \end{equation}

    \item Solve the Poisson equation to obtain pressure correction term
        \begin{equation}
          \begin{split}
            \frac{\partial}{\partial x^{k}}\frac{\partial p}{\partial x^{k}} = \frac{\partial u^{k}}{\partial x^{k}},\,p(x, 1) = 0,\,\frac{\partial p(x, 0)}{\partial y} = \frac{\partial p(0, y)}{\partial x} = \frac{\partial p(1, y)}{\partial x} = 0.
          \end{split}
        \end{equation}

    \item Correct speed to restore incompressibility
        \begin{equation}
            v^{i}_{(n+1)} = u^{i} - \frac{\partial p}{\partial x^{i}}
        \end{equation}
\end{enumerate}
Since the problem is defined on a deformed mesh, we need to rewrite the Chorin projection method in the curvilinear coordinates.

Under the coordinate transformations, the equation for pressure changes to
\begin{equation}
  \begin{split}
    &\frac{\partial}{\partial \xi^{i}}\left(J\delta^{\alpha\beta}\frac{\partial \xi^{i}}{\partial x^{\alpha}}\frac{\partial \xi^{j}}{\partial x^{\beta}}\frac{\partial p}{\partial \xi^{j}}\right) = \frac{\partial}{\partial \xi^{j}}\left(\overline{u}^{j} J\right),\,\overline{u}^{j} = \frac{\partial \xi^{j}}{\partial x^{\alpha}} u^{\alpha},\\
    & p\left(\xi^{1}, 1\right) = 1,\,\left(\frac{\partial p}{\partial \xi^{1}}\frac{\partial \xi^{1}}{\partial x^2} + \frac{\partial p}{\partial \xi^{2}}\frac{\partial \xi^{2}}{\partial x^2}\right)_{\xi^{2}=0} = 0,\\
    &\left(\frac{\partial p}{\partial \xi^{1}}\frac{\partial \xi^{1}}{\partial x^1} + \frac{\partial p}{\partial \xi^{2}}\frac{\partial \xi^{2}}{\partial x^1}\right)_{\xi^{1}=0} = 0,\,\left(\frac{\partial p}{\partial \xi^{1}}\frac{\partial \xi^{1}}{\partial x^1} + \frac{\partial p}{\partial \xi^{2}}\frac{\partial \xi^{2}}{\partial x^1}\right)_{\xi^{1}=1} = 0.
  \end{split}
\end{equation}

The equation for speed update is a vector conservation law
\begin{equation}
    \frac{\partial A^{ij}}{\partial x^{j}} = F^{i},
\end{equation}
with the following transformation rule (see \cite{MR3675718}, Section 2.4.2)
\begin{equation}
    \label{eq:vector conservation transform}
    \frac{\partial}{\partial \xi^{j}}\left(J \overline{A}^{ij}\right) + \frac{\partial^2 x^{l}}{\partial \xi^{k}\partial \xi^{j}}\frac{\partial \xi^{i}}{\partial x^{l}}\overline{A}^{kj} = \overline{F}^{i},\,\overline{A}^{kj} = \frac{\partial \xi^{k}}{\partial x^{\alpha}}\frac{\partial \xi^{j}}{\partial x^{\beta}}A^{\alpha\beta},\,\overline{F}^{j} = \frac{\partial \xi^{j}}{\partial x^{\alpha}}F^{\alpha}.
\end{equation}

It is straightforward to apply general equation \cref{eq:vector conservation transform} to a particular case
\begin{equation}
  \frac{u^{i} - v^{i}_{(n)}}{\Delta t} = \frac{\partial}{\partial x^{k}}\left(-v^{k}_{(n)} v^{i}_{(n)} + \nu \frac{\partial v^{i}_{(n)}}{\partial x^k}\right).
\end{equation}
The only problematic term is 
\begin{equation}
  \frac{\partial v^{i}_{(n)}}{\partial x^k} \longrightarrow \frac{\partial \xi^{i}}{\partial x^{\alpha}}\frac{\partial \xi^{k}}{\partial x^{\beta}} \frac{\partial v^{\alpha}_{(n)}}{\partial x^\beta} = \frac{\partial \xi^{i}}{\partial x^{\alpha}}\frac{\partial \xi^{k}}{\partial x^{\beta}}\frac{\partial \xi^{\rho}}{\partial x^{\beta}} \frac{\partial v^{\alpha}_{(n)}}{\partial \xi^\rho}.
\end{equation}
To evaluate this term, we need to switch from $\overline{v}^{\alpha}$ to $v^{\alpha}$ which is a simple task since the Jacobi matrix is available.
\section{Supplementray results}
\label{appendix:Supplementray results}
Here we present additional data on the experiments described in \cref{section:Experiments}.

\cref{table:elliptic_1D}, \cref{table:convection_diffusion_1D}, \cref{table:wave_10_1D}, \cref{table:wave_5_1D} contain results on $D=1$ experiments.

In \cref{table:sensitivity}, one can find data with sensitivity to coordinate transforms.

Results for DeepONet training on Navier-Stokes dataset are in \cref{table:DeepONet}.

\begin{table}[tbh]
    \centering
    \caption{Sensitivity to grid distortion for DilResNet and FNO with $\surd$ and without $\times$ augmentation. The distortion here refers to the maximal difference between the unperturbed $\boldsymbol{x}$ and  perturbed $\boldsymbol{x}(\boldsymbol{\xi})$ grids averaged over $1000$ grids used to augment dataset.}
    \vskip 0.15in
    \begin{tabular}{@{}lcccccccccccc@{}}
        \toprule
        & \multicolumn{4}{c}{Elliptic alpha} & \multicolumn{4}{c}{Convection-diffusion} & \multicolumn{4}{c}{Wave}\\\cmidrule(r){2-5}\cmidrule(r){6-9}\cmidrule(r){10-13}
        $\Delta$ & \multicolumn{2}{c}{DilResNet} & \multicolumn{2}{c}{FNO} & \multicolumn{2}{c}{DilResNet} & \multicolumn{2}{c}{FNO} & \multicolumn{2}{c}{DilResNet} & \multicolumn{2}{c}{FNO} \\ \addlinespace[0.1em]
         & $\times$ & $\surd$ & $\times$ & $\surd$ & $\times$ & $\surd$ & $\times$ & $\surd$ & $\times$ & $\surd$ & $\times$ & $\surd$ \\\addlinespace[0.5em]
        $0.006$ & $10\%$ & $2\%$ & $4\%$ & $1\%$ & $2\%$ & $1\%$ & $2\%$ & $1\%$ & $2\%$ & $2\%$ & $5\%$ & $2\%$\\
        $0.040$ & $14\%$ & $2\%$ & $7\%$ & $2\%$ & $15\%$ & $2\%$ & $15\%$ & $4\%$ & $17\%$ & $4\%$ & $28\%$ & $9\%$\\
        $0.099$ & $42\%$ & $3\%$ & $15\%$ & $5\%$ & --- & $25\%$ & --- & $26\%$ & --- & $18\%$ & --- & $34\%$\\
        $0.119$ & $70\%$ & $4\%$ & $18\%$ & $6\%$ & --- & $58\%$ & --- & $54\%$ & --- & $44\%$ & --- & $69\%$\\
        \bottomrule
    \end{tabular}
    \label{table:sensitivity}
    \vskip -0.1in
\end{table}
\begin{table}[ht]
\centering
\caption{Average test errors $\pm$ standard deviation for elliptic equation in one dimension. Factor $m$ in columns corresponds to the number of extra samples $m \times N_{\text{train}}$ added to the dataset with augmentation or resampling.}
\vskip 0.15in
\resizebox{1\textwidth}{!}{\begin{tabular}{@{}lcccccc@{}}
    \toprule
    Model&&$N_{\text{train}}\backslash m$ & $1$ & $2$ & $3$ & $4$ \\
    \cmidrule(r){4-7}
	\multirow{8}{*}{DeepONet}& \multirow{4}{*}{augmentation}&$500$&$0.573\pm0.145$&$0.419\pm0.032$&$0.415\pm0.038$&$0.427\pm0.041$\\
&&$1000$&$0.489\pm0.071$&$0.415\pm0.033$&$0.383\pm0.022$&$0.373\pm0.013$\\
&&$1500$&$0.461\pm0.041$&$0.398\pm0.019$&$0.399\pm0.075$&$0.465\pm0.101$\\
&&$2000$&$0.417\pm0.023$&$0.385\pm0.025$&$0.383\pm0.023$&$0.375\pm0.024$\\\cmidrule(r){4-7}
	 & \multirow{4}{*}{resampling}&$500$&$0.633\pm0.066$&$0.492\pm0.042$&$0.445\pm0.026$&$0.388\pm0.017$\\
&&$1000$&$0.597\pm0.07$&$0.536\pm0.072$&$0.441\pm0.023$&$0.705\pm0.606$\\
&&$1500$&$0.646\pm0.092$&$0.505\pm0.031$&$0.46\pm0.039$&$0.415\pm0.011$\\
&&$2000$&$0.624\pm0.072$&$0.504\pm0.04$&$0.428\pm0.026$&$0.452\pm0.078$\\\cmidrule(r){2-7}
	\multirow{8}{*}{FNO}& \multirow{4}{*}{augmentation}&$500$&$0.125\pm0.014$&$0.063\pm0.003$&$0.047\pm0.002$&$0.039\pm0.002$\\
&&$1000$&$0.09\pm0.006$&$0.052\pm0.002$&$0.04\pm0.002$&$0.032\pm0.001$\\
&&$1500$&$0.074\pm0.004$&$0.043\pm0.003$&$0.035\pm0.002$&$0.028\pm0.001$\\
&&$2000$&$0.064\pm0.004$&$0.04\pm0.002$&$0.032\pm0.002$&$0.026\pm0.001$\\\cmidrule(r){4-7}
	 & \multirow{4}{*}{resampling}&$500$&$0.12\pm0.001$&$0.068\pm0.004$&$0.053\pm0.002$&$0.043\pm0.002$\\
&&$1000$&$0.106\pm0.004$&$0.063\pm0.001$&$0.049\pm0.002$&$0.04\pm0.001$\\
&&$1500$&$0.102\pm0.005$&$0.061\pm0.003$&$0.049\pm0.002$&$0.041\pm0.003$\\
&&$2000$&$0.098\pm0.004$&$0.063\pm0.004$&$0.048\pm0.003$&$0.04\pm0.001$\\\cmidrule(r){2-7}
	\multirow{8}{*}{rFNO}& \multirow{4}{*}{augmentation}&$500$&$0.146\pm0.004$&$0.121\pm0.004$&$0.106\pm0.004$&$0.099\pm0.003$\\
&&$1000$&$0.103\pm0.002$&$0.087\pm0.004$&$0.08\pm0.004$&$0.076\pm0.002$\\
&&$1500$&$0.082\pm0.002$&$0.074\pm0.002$&$0.07\pm0.003$&$0.064\pm0.001$\\
&&$2000$&$0.073\pm0.001$&$0.065\pm0.002$&$0.061\pm0.001$&$0.056\pm0.001$\\\cmidrule(r){4-7}
	 & \multirow{4}{*}{resampling}&$500$&$0.17\pm0.005$&$0.154\pm0.003$&$0.148\pm0.006$&$0.148\pm0.004$\\
&&$1000$&$0.111\pm0.002$&$0.111\pm0.002$&$0.107\pm0.002$&$0.105\pm0.004$\\
&&$1500$&$0.089\pm0.002$&$0.087\pm0.003$&$0.086\pm0.002$&$0.084\pm0.002$\\
&&$2000$&$0.078\pm0.002$&$0.076\pm0.003$&$0.074\pm0.002$&$0.075\pm0.002$\\\cmidrule(r){2-7}
	\multirow{8}{*}{DilResNet}& \multirow{4}{*}{augmentation}&$500$&$0.374\pm0.037$&$0.304\pm0.033$&$0.243\pm0.031$&$0.214\pm0.026$\\
&&$1000$&$0.208\pm0.005$&$0.172\pm0.013$&$0.145\pm0.004$&$0.141\pm0.015$\\
&&$1500$&$0.179\pm0.013$&$0.143\pm0.006$&$0.119\pm0.009$&$0.114\pm0.01$\\
&&$2000$&$0.13\pm0.009$&$0.118\pm0.009$&$0.104\pm0.008$&$0.095\pm0.007$\\\cmidrule(r){4-7}
	 & \multirow{4}{*}{resampling}&$500$&$0.485\pm0.044$&$0.48\pm0.087$&$0.41\pm0.063$&$0.425\pm0.052$\\
&&$1000$&$0.255\pm0.009$&$0.255\pm0.027$&$0.218\pm0.011$&$0.226\pm0.021$\\
&&$1500$&$0.185\pm0.024$&$0.17\pm0.014$&$0.171\pm0.008$&$0.17\pm0.015$\\
&&$2000$&$0.15\pm0.014$&$0.141\pm0.011$&$0.143\pm0.008$&$0.14\pm0.005$\\\cmidrule(r){2-7}
	\multirow{8}{*}{MLP}& \multirow{4}{*}{augmentation}&$500$&$0.346\pm0.031$&$0.337\pm0.031$&$0.261\pm0.041$&$0.262\pm0.043$\\
&&$1000$&$0.268\pm0.042$&$0.181\pm0.03$&$0.162\pm0.026$&$0.126\pm0.016$\\
&&$1500$&$0.23\pm0.049$&$0.145\pm0.017$&$0.13\pm0.007$&$0.105\pm0.009$\\
&&$2000$&$0.155\pm0.028$&$0.109\pm0.007$&$0.09\pm0.005$&$0.089\pm0.012$\\\cmidrule(r){4-7}
	 & \multirow{4}{*}{resampling}&$500$&$0.342\pm0.023$&$0.377\pm0.011$&$0.332\pm0.066$&$0.286\pm0.053$\\
&&$1000$&$0.309\pm0.048$&$0.301\pm0.105$&$0.239\pm0.096$&$0.184\pm0.043$\\
&&$1500$&$0.211\pm0.03$&$0.133\pm0.026$&$0.111\pm0.007$&$0.111\pm0.018$\\
&&$2000$&$0.151\pm0.038$&$0.099\pm0.013$&$0.096\pm0.012$&$0.089\pm0.006$\\
\bottomrule
\end{tabular}}
\label{table:elliptic_1D}
\vskip -0.1in
\end{table}
\begin{table}[ht]
\centering
\caption{Average test errors $\pm$ standard deviation for convection-diffusion equation in one dimension. Factor $m$ in columns corresponds to the number of extra samples $m\times N_{\text{train}}$ added to the dataset with augmentation or resampling.}
\vskip 0.15in
\resizebox{1\textwidth}{!}{\begin{tabular}{@{}lcccccc@{}}
    \toprule
    Model&&$N_{\text{train}}\backslash m$ & $1$ & $2$ & $3$ & $4$ \\
    \cmidrule(r){4-7}
	\multirow{8}{*}{DeepONet}& \multirow{4}{*}{augmentation}&$500$&$0.8\pm0.013$&$0.619\pm0.012$&$0.522\pm0.022$&$0.449\pm0.018$\\
&&$1000$&$0.716\pm0.013$&$0.536\pm0.004$&$0.469\pm0.013$&$0.401\pm0.02$\\
&&$1500$&$0.646\pm0.02$&$0.477\pm0.009$&$0.407\pm0.01$&$0.327\pm0.02$\\
&&$2000$&$0.607\pm0.019$&$0.46\pm0.014$&$0.386\pm0.015$&$0.312\pm0.019$\\\cmidrule(r){4-7}
	 & \multirow{4}{*}{resampling}&$500$&$0.92\pm0.034$&$0.771\pm0.02$&$0.671\pm0.025$&$0.619\pm0.024$\\
&&$1000$&$0.915\pm0.017$&$0.769\pm0.016$&$0.665\pm0.017$&$0.612\pm0.022$\\
&&$1500$&$0.926\pm0.041$&$0.772\pm0.017$&$0.676\pm0.011$&$0.619\pm0.02$\\
&&$2000$&$0.92\pm0.041$&$0.768\pm0.031$&$0.672\pm0.019$&$0.623\pm0.03$\\\cmidrule(r){2-7}
	\multirow{8}{*}{FNO}& \multirow{4}{*}{augmentation}&$500$&$0.464\pm0.008$&$0.366\pm0.023$&$0.235\pm0.039$&$0.117\pm0.026$\\
&&$1000$&$0.472\pm0.011$&$0.324\pm0.033$&$0.137\pm0.036$&$0.072\pm0.007$\\
&&$1500$&$0.438\pm0.043$&$0.229\pm0.085$&$0.101\pm0.017$&$0.059\pm0.005$\\
&&$2000$&$0.426\pm0.036$&$0.212\pm0.057$&$0.081\pm0.016$&$0.05\pm0.008$\\\cmidrule(r){4-7}
	 & \multirow{4}{*}{resampling}&$500$&$0.499\pm0.01$&$0.426\pm0.016$&$0.315\pm0.027$&$0.157\pm0.031$\\
&&$1000$&$0.529\pm0.009$&$0.442\pm0.026$&$0.258\pm0.04$&$0.133\pm0.028$\\
&&$1500$&$0.553\pm0.013$&$0.434\pm0.03$&$0.264\pm0.051$&$0.124\pm0.017$\\
&&$2000$&$0.543\pm0.009$&$0.444\pm0.02$&$0.233\pm0.033$&$0.125\pm0.013$\\\cmidrule(r){2-7}
	\multirow{8}{*}{rFNO}& \multirow{4}{*}{augmentation}&$500$&$0.524\pm0.008$&$0.49\pm0.007$&$0.436\pm0.039$&$0.402\pm0.049$\\
&&$1000$&$0.176\pm0.005$&$0.157\pm0.037$&$0.124\pm0.006$&$0.115\pm0.007$\\
&&$1500$&$0.108\pm0.003$&$0.097\pm0.004$&$0.088\pm0.004$&$0.082\pm0.004$\\
&&$2000$&$0.083\pm0.005$&$0.071\pm0.002$&$0.068\pm0.002$&$0.067\pm0.004$\\\cmidrule(r){4-7}
	 & \multirow{4}{*}{resampling}&$500$&$0.536\pm0.007$&$0.513\pm0.006$&$0.507\pm0.007$&$0.485\pm0.019$\\
&&$1000$&$0.265\pm0.057$&$0.258\pm0.04$&$0.219\pm0.013$&$0.194\pm0.02$\\
&&$1500$&$0.137\pm0.008$&$0.124\pm0.009$&$0.119\pm0.005$&$0.115\pm0.007$\\
&&$2000$&$0.102\pm0.004$&$0.094\pm0.003$&$0.09\pm0.003$&$0.09\pm0.002$\\\cmidrule(r){2-7}
	\multirow{8}{*}{DilResNet}& \multirow{4}{*}{augmentation}&$500$&$0.133\pm0.014$&$0.109\pm0.006$&$0.1\pm0.006$&$0.091\pm0.008$\\
&&$1000$&$0.07\pm0.004$&$0.058\pm0.002$&$0.048\pm0.004$&$0.044\pm0.002$\\
&&$1500$&$0.047\pm0.004$&$0.041\pm0.002$&$0.037\pm0.002$&$0.033\pm0.003$\\
&&$2000$&$0.038\pm0.003$&$0.032\pm0.001$&$0.029\pm0.001$&$0.026\pm0.002$\\\cmidrule(r){4-7}
	 & \multirow{4}{*}{resampling}&$500$&$0.155\pm0.006$&$0.172\pm0.026$&$0.144\pm0.008$&$0.144\pm0.014$\\
&&$1000$&$0.075\pm0.002$&$0.072\pm0.005$&$0.067\pm0.003$&$0.068\pm0.007$\\
&&$1500$&$0.052\pm0.002$&$0.05\pm0.003$&$0.047\pm0.003$&$0.047\pm0.001$\\
&&$2000$&$0.041\pm0.002$&$0.038\pm0.002$&$0.036\pm0.002$&$0.037\pm0.001$\\\cmidrule(r){2-7}
	\multirow{8}{*}{MLP}& \multirow{4}{*}{augmentation}&$500$&$0.425\pm0.03$&$0.439\pm0.011$&$0.439\pm0.031$&$0.423\pm0.015$\\
&&$1000$&$0.306\pm0.049$&$0.301\pm0.045$&$0.291\pm0.086$&$0.266\pm0.055$\\
&&$1500$&$0.26\pm0.059$&$0.175\pm0.058$&$0.125\pm0.013$&$0.115\pm0.032$\\
&&$2000$&$0.105\pm0.008$&$0.101\pm0.015$&$0.071\pm0.005$&$0.074\pm0.014$\\\cmidrule(r){4-7}
	 & \multirow{4}{*}{resampling}&$500$&$0.464\pm0.018$&$0.507\pm0.037$&$0.5\pm0.014$&$0.482\pm0.012$\\
&&$1000$&$0.381\pm0.034$&$0.361\pm0.086$&$0.379\pm0.047$&$0.328\pm0.061$\\
&&$1500$&$0.306\pm0.056$&$0.299\pm0.064$&$0.2\pm0.074$&$0.231\pm0.039$\\
&&$2000$&$0.19\pm0.052$&$0.183\pm0.091$&$0.135\pm0.023$&$0.114\pm0.015$\\
\bottomrule
\end{tabular}}
\label{table:convection_diffusion_1D}
\vskip -0.1in
\end{table}
\begin{table}[ht]
\centering
\caption{Average test errors $\pm$ standard deviation for wave equation (10 modes) in one dimension. Factor $m$ in columns corresponds to the number of extra samples $m\times N_{\text{train}}$ added to the dataset with augmentation or resampling.}
\vskip 0.15in
\resizebox{1\textwidth}{!}{\begin{tabular}{@{}lcccccc@{}}
    \toprule
    Model&&$N_{\text{train}}\backslash m$ & $1$ & $2$ & $3$ & $4$ \\
    \cmidrule(r){4-7}
	\multirow{8}{*}{DeepONet}& \multirow{4}{*}{augmentation}&$500$&$0.275\pm0.009$&$0.215\pm0.004$&$0.181\pm0.005$&$0.171\pm0.005$\\
&&$1000$&$0.259\pm0.004$&$0.207\pm0.005$&$0.175\pm0.004$&$0.163\pm0.004$\\
&&$1500$&$0.255\pm0.001$&$0.202\pm0.003$&$0.174\pm0.007$&$0.162\pm0.006$\\
&&$2000$&$0.244\pm0.003$&$0.197\pm0.005$&$0.172\pm0.003$&$0.164\pm0.006$\\\cmidrule(r){4-7}
	 & \multirow{4}{*}{resampling}&$500$&$0.322\pm0.014$&$0.241\pm0.004$&$0.198\pm0.003$&$0.18\pm0.006$\\
&&$1000$&$0.319\pm0.009$&$0.25\pm0.021$&$0.196\pm0.005$&$0.179\pm0.005$\\
&&$1500$&$0.314\pm0.011$&$0.237\pm0.005$&$0.212\pm0.029$&$0.18\pm0.007$\\
&&$2000$&$0.321\pm0.008$&$0.237\pm0.005$&$0.197\pm0.005$&$0.18\pm0.006$\\\cmidrule(r){2-7}
	\multirow{8}{*}{FNO}& \multirow{4}{*}{augmentation}&$500$&$0.135\pm0.003$&$0.088\pm0.002$&$0.068\pm0.002$&$0.057\pm0.001$\\
&&$1000$&$0.123\pm0.003$&$0.082\pm0.003$&$0.062\pm0.001$&$0.049\pm0.001$\\
&&$1500$&$0.119\pm0.004$&$0.076\pm0.002$&$0.058\pm0.001$&$0.045\pm0.001$\\
&&$2000$&$0.108\pm0.002$&$0.072\pm0.002$&$0.052\pm0.001$&$0.042\pm0.001$\\\cmidrule(r){4-7}
	 & \multirow{4}{*}{resampling}&$500$&$0.208\pm0.008$&$0.14\pm0.007$&$0.108\pm0.004$&$0.091\pm0.003$\\
&&$1000$&$0.178\pm0.004$&$0.124\pm0.004$&$0.096\pm0.004$&$0.078\pm0.005$\\
&&$1500$&$0.168\pm0.003$&$0.118\pm0.005$&$0.09\pm0.004$&$0.074\pm0.003$\\
&&$2000$&$0.163\pm0.003$&$0.114\pm0.003$&$0.087\pm0.003$&$0.07\pm0.002$\\\cmidrule(r){2-7}
	\multirow{8}{*}{rFNO}& \multirow{4}{*}{augmentation}&$500$&$0.19\pm0.003$&$0.179\pm0.002$&$0.172\pm0.003$&$0.166\pm0.003$\\
&&$1000$&$0.147\pm0.002$&$0.137\pm0.002$&$0.131\pm0.002$&$0.127\pm0.003$\\
&&$1500$&$0.127\pm0.003$&$0.117\pm0.004$&$0.112\pm0.002$&$0.108\pm0.002$\\
&&$2000$&$0.111\pm0.001$&$0.103\pm0.001$&$0.098\pm0.002$&$0.094\pm0.002$\\\cmidrule(r){4-7}
	 & \multirow{4}{*}{resampling}&$500$&$0.213\pm0.004$&$0.212\pm0.003$&$0.21\pm0.002$&$0.207\pm0.005$\\
&&$1000$&$0.168\pm0.005$&$0.168\pm0.002$&$0.164\pm0.003$&$0.163\pm0.002$\\
&&$1500$&$0.144\pm0.004$&$0.141\pm0.001$&$0.142\pm0.003$&$0.143\pm0.003$\\
&&$2000$&$0.129\pm0.003$&$0.127\pm0.003$&$0.126\pm0.002$&$0.126\pm0.002$\\\cmidrule(r){2-7}
	\multirow{8}{*}{DilResNet}& \multirow{4}{*}{augmentation}&$500$&$0.157\pm0.009$&$0.133\pm0.007$&$0.128\pm0.006$&$0.121\pm0.005$\\
&&$1000$&$0.109\pm0.006$&$0.092\pm0.003$&$0.086\pm0.007$&$0.084\pm0.005$\\
&&$1500$&$0.084\pm0.006$&$0.076\pm0.003$&$0.072\pm0.004$&$0.066\pm0.005$\\
&&$2000$&$0.076\pm0.007$&$0.063\pm0.004$&$0.058\pm0.004$&$0.056\pm0.005$\\\cmidrule(r){4-7}
	 & \multirow{4}{*}{resampling}&$500$&$0.173\pm0.008$&$0.174\pm0.007$&$0.175\pm0.012$&$0.182\pm0.013$\\
&&$1000$&$0.121\pm0.007$&$0.119\pm0.011$&$0.118\pm0.009$&$0.126\pm0.014$\\
&&$1500$&$0.098\pm0.006$&$0.095\pm0.006$&$0.098\pm0.007$&$0.093\pm0.011$\\
&&$2000$&$0.082\pm0.005$&$0.08\pm0.002$&$0.08\pm0.005$&$0.079\pm0.003$\\\cmidrule(r){2-7}
	\multirow{8}{*}{MLP}& \multirow{4}{*}{augmentation}&$500$&$0.349\pm0.038$&$0.331\pm0.041$&$0.304\pm0.041$&$0.294\pm0.015$\\
&&$1000$&$0.231\pm0.022$&$0.198\pm0.021$&$0.178\pm0.015$&$0.187\pm0.033$\\
&&$1500$&$0.186\pm0.013$&$0.145\pm0.009$&$0.137\pm0.007$&$0.134\pm0.013$\\
&&$2000$&$0.151\pm0.012$&$0.131\pm0.01$&$0.115\pm0.012$&$0.099\pm0.006$\\\cmidrule(r){4-7}
	 & \multirow{4}{*}{resampling}&$500$&$0.376\pm0.072$&$0.349\pm0.044$&$0.341\pm0.033$&$0.334\pm0.052$\\
&&$1000$&$0.236\pm0.01$&$0.225\pm0.015$&$0.236\pm0.036$&$0.233\pm0.017$\\
&&$1500$&$0.188\pm0.005$&$0.192\pm0.012$&$0.178\pm0.009$&$0.193\pm0.01$\\
&&$2000$&$0.152\pm0.007$&$0.159\pm0.004$&$0.158\pm0.001$&$0.17\pm0.014$\\
\bottomrule
\end{tabular}}
\label{table:wave_10_1D}
\vskip -0.1in
\end{table}
\begin{table}[ht]
\centering
\caption{Average test errors $\pm$ standard deviation for wave equation (5 modes) in one dimension. Factor $m$ in columns corresponds to the number of extra samples $m\times N_{\text{train}}$ added to the dataset with augmentation or resampling.}
\vskip 0.15in
\resizebox{1\textwidth}{!}{\begin{tabular}{@{}lcccccc@{}}
    \toprule
    Model&&$N_{\text{train}}\backslash m$ & $1$ & $2$ & $3$ & $4$ \\
    \cmidrule(r){4-7}
	\multirow{8}{*}{DeepONet}& \multirow{4}{*}{augmentation}&$500$&$0.262\pm0.005$&$0.178\pm0.004$&$0.152\pm0.003$&$0.147\pm0.011$\\
&&$1000$&$0.247\pm0.011$&$0.169\pm0.006$&$0.158\pm0.013$&$0.141\pm0.008$\\
&&$1500$&$0.237\pm0.006$&$0.161\pm0.008$&$0.143\pm0.006$&$0.134\pm0.004$\\
&&$2000$&$0.228\pm0.005$&$0.162\pm0.004$&$0.147\pm0.01$&$0.138\pm0.005$\\\cmidrule(r){4-7}
	 & \multirow{4}{*}{resampling}&$500$&$0.306\pm0.008$&$0.196\pm0.004$&$0.157\pm0.005$&$0.145\pm0.008$\\
&&$1000$&$0.31\pm0.004$&$0.195\pm0.004$&$0.159\pm0.006$&$0.148\pm0.008$\\
&&$1500$&$0.308\pm0.011$&$0.197\pm0.004$&$0.172\pm0.015$&$0.144\pm0.006$\\
&&$2000$&$0.303\pm0.006$&$0.197\pm0.007$&$0.16\pm0.007$&$0.142\pm0.005$\\\cmidrule(r){2-7}
	\multirow{8}{*}{FNO}& \multirow{4}{*}{augmentation}&$500$&$0.088\pm0.002$&$0.057\pm0.002$&$0.042\pm0.002$&$0.034\pm0.001$\\
&&$1000$&$0.083\pm0.003$&$0.051\pm0.002$&$0.038\pm0.001$&$0.031\pm0.001$\\
&&$1500$&$0.078\pm0.002$&$0.048\pm0.001$&$0.036\pm0.001$&$0.029\pm0.001$\\
&&$2000$&$0.074\pm0.0$&$0.045\pm0.001$&$0.033\pm0.001$&$0.027\pm0.001$\\\cmidrule(r){4-7}
	 & \multirow{4}{*}{resampling}&$500$&$0.18\pm0.012$&$0.104\pm0.009$&$0.073\pm0.005$&$0.058\pm0.002$\\
&&$1000$&$0.138\pm0.007$&$0.081\pm0.004$&$0.057\pm0.003$&$0.047\pm0.002$\\
&&$1500$&$0.12\pm0.004$&$0.071\pm0.003$&$0.052\pm0.002$&$0.041\pm0.002$\\
&&$2000$&$0.109\pm0.004$&$0.065\pm0.003$&$0.047\pm0.001$&$0.039\pm0.001$\\\cmidrule(r){2-7}
	\multirow{8}{*}{rFNO}& \multirow{4}{*}{augmentation}&$500$&$0.16\pm0.003$&$0.145\pm0.003$&$0.138\pm0.002$&$0.136\pm0.002$\\
&&$1000$&$0.11\pm0.003$&$0.102\pm0.002$&$0.098\pm0.002$&$0.094\pm0.002$\\
&&$1500$&$0.092\pm0.002$&$0.085\pm0.002$&$0.08\pm0.001$&$0.077\pm0.002$\\
&&$2000$&$0.078\pm0.002$&$0.072\pm0.001$&$0.069\pm0.001$&$0.065\pm0.001$\\\cmidrule(r){4-7}
	 & \multirow{4}{*}{resampling}&$500$&$0.176\pm0.004$&$0.175\pm0.004$&$0.174\pm0.002$&$0.173\pm0.002$\\
&&$1000$&$0.125\pm0.005$&$0.126\pm0.001$&$0.122\pm0.002$&$0.123\pm0.002$\\
&&$1500$&$0.104\pm0.003$&$0.101\pm0.002$&$0.1\pm0.001$&$0.1\pm0.001$\\
&&$2000$&$0.089\pm0.001$&$0.087\pm0.001$&$0.087\pm0.002$&$0.087\pm0.002$\\\cmidrule(r){2-7}
	\multirow{8}{*}{DilResNet}& \multirow{4}{*}{augmentation}&$500$&$0.132\pm0.012$&$0.113\pm0.01$&$0.116\pm0.014$&$0.105\pm0.009$\\
&&$1000$&$0.093\pm0.008$&$0.091\pm0.009$&$0.068\pm0.005$&$0.068\pm0.003$\\
&&$1500$&$0.077\pm0.008$&$0.066\pm0.009$&$0.061\pm0.006$&$0.055\pm0.006$\\
&&$2000$&$0.061\pm0.007$&$0.053\pm0.004$&$0.05\pm0.003$&$0.055\pm0.004$\\\cmidrule(r){4-7}
	 & \multirow{4}{*}{resampling}&$500$&$0.169\pm0.019$&$0.174\pm0.013$&$0.147\pm0.009$&$0.155\pm0.009$\\
&&$1000$&$0.098\pm0.007$&$0.105\pm0.01$&$0.112\pm0.008$&$0.11\pm0.009$\\
&&$1500$&$0.09\pm0.011$&$0.082\pm0.005$&$0.076\pm0.003$&$0.079\pm0.005$\\
&&$2000$&$0.065\pm0.003$&$0.07\pm0.006$&$0.068\pm0.006$&$0.066\pm0.006$\\\cmidrule(r){2-7}
	\multirow{8}{*}{MLP}& \multirow{4}{*}{augmentation}&$500$&$0.395\pm0.09$&$0.277\pm0.013$&$0.286\pm0.016$&$0.243\pm0.015$\\
&&$1000$&$0.259\pm0.049$&$0.169\pm0.025$&$0.144\pm0.022$&$0.131\pm0.017$\\
&&$1500$&$0.172\pm0.02$&$0.138\pm0.034$&$0.098\pm0.009$&$0.085\pm0.005$\\
&&$2000$&$0.112\pm0.014$&$0.084\pm0.007$&$0.076\pm0.005$&$0.069\pm0.007$\\\cmidrule(r){4-7}
	 & \multirow{4}{*}{resampling}&$500$&$0.358\pm0.051$&$0.387\pm0.053$&$0.383\pm0.057$&$0.401\pm0.054$\\
&&$1000$&$0.288\pm0.122$&$0.224\pm0.018$&$0.202\pm0.029$&$0.233\pm0.049$\\
&&$1500$&$0.178\pm0.011$&$0.155\pm0.015$&$0.14\pm0.031$&$0.136\pm0.011$\\
&&$2000$&$0.148\pm0.032$&$0.119\pm0.018$&$0.103\pm0.01$&$0.111\pm0.013$\\
\bottomrule
\end{tabular}}
\label{table:wave_5_1D}
\vskip -0.1in
\end{table}
\begin{table}[tbh]
    \centering
    \caption{Relative errors for DeepONet, Navier-Stokes dataset, $\surd$ marks training run with augmentation and $\times$ --- without augmentation.}
    \vskip 0.15in
    \begin{tabular}{@{}lcccccccc@{}}
        \toprule
        & \multicolumn{4}{c}{$v^{1}$} & \multicolumn{4}{c}{$v^{2}$}\\\cmidrule(r){2-5}\cmidrule(r){6-9}
        model & \multicolumn{2}{c}{$E_{\text{train}}$} & \multicolumn{2}{c}{$E_{\text{test}}$} & \multicolumn{2}{c}{$E_{\text{train}}$} & \multicolumn{2}{c}{$E_{\text{test}}$} \\ \addlinespace[0.1em]
         & $\times$ & $\surd$ & $\times$ & $\surd$ & $\times$ & $\surd$ & $\times$ & $\surd$ \\\addlinespace[0.5em]
        DeepONet & $0.074$ & $0.063$ & $0.162$ & $0.163$ & $0.212$ & $0.232$ & $0.368$ & $0.395$ \\
        POD-DeepONet & $0.622$ & $0.622$ & $0.589$ & $0.589$ & $0.349$ & $0.358$ & $0.414$ & $0.419$ \\
        \bottomrule
    \end{tabular}
    \label{table:DeepONet}
    \vskip -0.1in
\end{table}

\end{document}